\documentclass[11pt]{article}

\usepackage[margin=1in]{geometry}
\usepackage{amsmath,amssymb,amsthm,bm}
\usepackage{booktabs}
\usepackage{graphicx}
\usepackage{placeins}
\usepackage{tikz}
\usetikzlibrary{arrows.meta,positioning,calc}
\graphicspath{{figures/}{../figures/}}
\usepackage[numbers,sort&compress]{natbib}
\usepackage[colorlinks=true,citecolor=blue,linkcolor=blue,urlcolor=blue]{hyperref}

\newtheorem{theorem}{Theorem}
\newtheorem{proposition}{Proposition}
\newtheorem{corollary}{Corollary}

\newcommand{\E}{\mathbb{E}}
\newcommand{\Prb}{\mathbb{P}}

\begin{document}

\title{Pooling Mobility Obscures Epidemic Invasion Routes}
\author{Tiandong Wang\thanks{Corresponding author: td\_wang@fudan.edu.cn}\\
Shanghai Center for Mathematical Sciences, Fudan University, Shanghai, China
\and
Wei Yang\\
Research Institute of Intelligent Complex Systems, Fudan University, Shanghai, China}

\maketitle

\begin{abstract}
Epidemic models often pool air travel, commuting, and other mobility layers
into a single weighted network. Pooling keeps the total imported infections into
a region but discards the transport mode and route that delivered them, the
information a mode-specific intervention needs. We make this precise for directed
multilayer flows with heavy-tailed variation. Pooling acts as an asymmetric
filter. The magnitude of an extreme importation, its radial tail index, is an
exact invariant of sampling and nonnegative aggregation, whereas its composition
across layers and routes, the angular part, is reweighted by layer-specific
sampling and can be irrecoverably merged. We give a necessary and sufficient
condition for recovering route composition, attach a surveillance decision cost
to the loss, and show that the first established route need not be the busiest. A controlled metapopulation
study calibrates the cost, and two contrasting reconstructions show where it
bites. In US pandemic influenza, air adds fitted information beyond commuting,
and a layer-resolved ranking targets more air-import burden than a pooled one.
In the early Italian COVID-19 wave, commuting improves onset reconstruction
while air attribution stays unresolved. Pooling can therefore support
total-importation surveillance but cannot in general identify the mode or route
behind an importation.
\end{abstract}

\noindent\textbf{Keywords:} multilayer networks; epidemic invasion; mobility;
regular variation; aggregation; surveillance

\section{Introduction}

Many complex systems transmit mass, information, or activity through several
types of directed connection. Multilayer network representations retain these
channels explicitly, whereas a common reduction replaces the layer tensor
by a single weighted adjacency matrix
\citep{dedomenico2013multilayer,kivela2014multilayer,boccaletti2014multilayer}.
We use \emph{coarse-graining} for this general reduction in resolution.
\emph{Aggregation} denotes the mathematical map that implements the reduction,
and \emph{pooling} denotes the specific operation of summing mobility layers or
route counts.
That reduction is often useful: it preserves total flow by construction and can
greatly simplify dynamics and inference. It can also identify network states
that are distinct before aggregation. We therefore ask which features are
mathematical invariants of the aggregation map and which become
nonidentifiable after aggregation.

The distinction is especially sharp for heterogeneous networks with extreme
flows. A multivariate extreme has a radial component, which measures total
magnitude, and an angular component, which records how that magnitude is
distributed across layers or routes. In the epidemic reading, the radial
component is the total infectious importation into a region, and the angular
component is the mix of transport modes and source regions that delivered it.
Summation can preserve radial tail order while mapping many angular directions to
the same aggregate direction. A second transformation acts downstream, where
network intensities become realized counts through thinning, random rescaling, or
mixed-Poisson sampling. We ask which radial and angular features survive this
sequence of stochastic transfer and deterministic aggregation.

We answer this question for directed multilayer networks with heavy-tailed
structural weights. Multilayer inhomogeneous random graphs and multivariate
regular variation provide the structural language
\citep{cirkovic2025mirg,han2025utd}, while generalized Breiman scaling and
Poissonization control the stochastic transfer
\citep{wang2024heterogeneous,resnick2007}. First, a common heavy-tailed factor
generates multivariate regular variation and an explicit interlayer upper-tail
dependence coefficient. Second, a random diagonal map followed by mixed-Poisson
sampling preserves the radial index but pushes forward the angular measure.
Third, route aggregation preserves the radial index but retains only a
norm-weighted pushforward of route composition. Whenever the normalized
aggregation map is noninjective, distinct route-level spectral measures are
observationally equivalent. Finally, competing route clocks connect the same
spectral measure to the first successful propagated event and its time scale.
We use the regular-variation formalism as an instrument. It makes
precise which surveillance-relevant quantities a pooled network still determines,
which it cannot recover, and what the loss costs a decision maker.

These results are sharper than the familiar statement that aggregation discards
multilayer structure. That structure is lost is well known, and
structural-reducibility methods choose how many layers to keep by minimizing an
information or redundancy criterion \citep{dedomenico2015reducibility}. For
extreme flows the loss is exact and asymmetric. The radial tail index is an exact
invariant of nonnegative aggregation, so pooling preserves the scaling of extreme
magnitude, whereas the angular spectral measure that encodes route and mode is
generically nonidentifiable after pooling. Proposition~\ref{prop:identifiability}
gives the necessary and sufficient condition under which route composition is
recoverable, and it identifies the pooled targets, such as total importation into
a destination, that stay exact. The loss carries a decision cost, the
surveillance regret of \eqref{eq:regret}, and a dynamical consequence through the
first-invasion clock, on which the first established route need not be the
busiest. Effective-distance approaches predict arrival times from aggregate
mobility magnitude \citep{brockmann2013}, using exactly the radial information
that pooling preserves. The route and mode attribution that pooling destroys is
the complementary quantity, and it is what the present results make precise.

Human mobility and epidemic invasion provide a concrete application in which
the lost angular information has a direct interpretation. Air and commuting
flows form distinct directed layers, and propagated events retain a source,
destination, and transport mode. Metapopulation models commonly pool these
flows into one coupling matrix
\citep{colizza2007,colizza2008,balcan2009,viboud2006}, even though mobility
choice and spatial resolution can change inferred origins and predictive
performance
\citep{klochkov2026origins,pullano2026resolution,charu2017influenza,kramer2020forecasting,klamser2024import}.
Joint node prominence can arise from a common
demographic factor, from population-residual gateway structure, or from
temporally synchronized fluctuations. These mechanisms need not agree. A region
may be prominent in both layers because it is large, because it has a specialized
transport role, or because a transient shock aligns its flows with those of
other regions. Aggregating the layers erases the distinction even when it
preserves total exposure.

We test this principle in finite networks, a stochastic metapopulation model,
and two contrasting case studies, US pandemic influenza and the early COVID-19
wave in Italy. The controlled synthetic study carries the quantitative claim,
because it fixes the data-generating process and isolates one mechanism at a
time through demographic alignment, gateway alignment, spatial placement, and
layer aggregation. The two applications apply the same hazard form in different
mobility regimes. In US pandemic influenza, air adds fitted information beyond
commuting, and pooling misdirects a mode-specific air-surveillance ranking. In
the early Italian wave, commuting improves onset reconstruction while air
attribution stays unresolved, the same nonidentifiability seen in observational
data. Across these settings, aggregate networks retain the scale of extreme
movement while failing to identify the route or mechanism that generated it. We
collect the caveats that qualify the two applications in the Discussion.

\section{Model Setup}
\label{sec:model}

\subsection{Directed multilayer network}

Let $V=\{1,\ldots,n\}$ index nodes, let
$\bm{x}_i\in\mathcal X\subset\mathbb R^2$ be the location of node $i$, and let
$N_i\in(0,\infty)$ be a positive node-scale variable. Let $\mathcal L$ be a
finite layer set and let time be indexed by $t\in\mathcal T$. In the mobility
application, nodes are regional populations, $N_i$ is resident population, and
$\mathcal L=\{a,c\}$ denotes aviation and commuting. For $i\ne j$, the
integer-valued random variable
$A_{ij}^{(\ell)}(t)\in\mathbb N_0$ counts directed movement units from $i$ to
$j$ in layer $\ell$ during period $t$; set $A_{ii}^{(\ell)}(t)=0$.
Then $\mathcal G_\ell(t)=(V,A^{(\ell)}(t))$ is a directed weighted
multigraph, with parallel trips represented by their multiplicity. Its outward
and inward strengths are, respectively,
\[
D_{i\ell}^{\mathrm{out}}(t)=\sum_{j\ne i}A_{ij}^{(\ell)}(t),\qquad
D_{i\ell}^{\mathrm{in}}(t)=\sum_{j\ne i}A_{ji}^{(\ell)}(t).
\]

The mathematical construction separates three maps: latent node structure to
directed layer flow, layer flow to propagated route counts, and route counts to
aggregated observations. The first two are stochastic; the third is a
nonnegative linear map. This separation lets us distinguish changes in radial
magnitude from loss of angular composition.

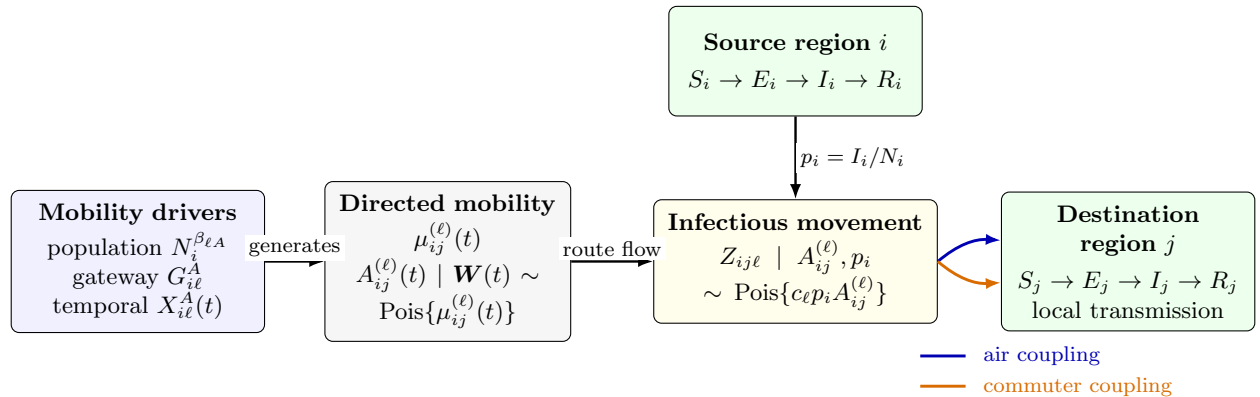
\begin{figure}[ht]
\centering
\begin{tikzpicture}[
  box/.style={draw, rounded corners=3pt, align=center, inner sep=5pt,
    minimum height=1.45cm, font=\footnotesize},
  flow/.style={-{Latex[length=2mm]}, line width=0.9pt},
  air/.style={-{Latex[length=2mm]}, line width=1.1pt, blue!70!black},
  commute/.style={-{Latex[length=2mm]}, line width=1.1pt, orange!85!black},
  note/.style={font=\scriptsize, fill=white, inner sep=1.5pt}
]
  \node[box, fill=blue!6, text width=3.0cm] (drivers) at (0,0) {
    \textbf{Mobility drivers}\\[2pt]
    population $N_i^{\beta_{\ell A}}$\\
    gateway $G_{i\ell}^{A}$\\
    temporal $X_{i\ell}^{A}(t)$};
  \node[box, fill=gray!7, text width=2.9cm] (mobility) at (4.1,0) {
    \textbf{Directed mobility}\\[2pt]
    $\mu_{ij}^{(\ell)}(t)$\\
    $A_{ij}^{(\ell)}(t)\mid\bm W(t)
      \sim\mathrm{Pois}\{\mu_{ij}^{(\ell)}(t)\}$};
  \node[box, fill=yellow!9, text width=3.4cm] (sampling) at (8.7,0) {
    \textbf{Infectious movement}\\[2pt]
    $Z_{ij\ell}\mid A_{ij}^{(\ell)},p_i$\\
    $\sim\mathrm{Pois}\{c_\ell p_iA_{ij}^{(\ell)}\}$};
  \node[box, fill=green!7, text width=3.0cm] (destination) at (13.1,0) {
    \textbf{Destination region $j$}\\[3pt]
    $S_j\rightarrow E_j\rightarrow I_j\rightarrow R_j$\\
    local transmission};
  \node[box, fill=green!7, text width=3.0cm] (source) at (8.7,2.65) {
    \textbf{Source region $i$}\\[3pt]
    $S_i\rightarrow E_i\rightarrow I_i\rightarrow R_i$};

  \draw[flow] (drivers.east) -- node[note, above] {generates} (mobility.west);
  \draw[flow] (mobility.east) -- node[note, above] {route flow} (sampling.west);
  \draw[flow] (source.south) -- node[note, right] {$p_i=I_i/N_i$} (sampling.north);
  \draw[air] (sampling.east) to[bend left=22]
    ([yshift=0.28cm]destination.west);
  \draw[commute] (sampling.east) to[bend right=22]
    ([yshift=-0.28cm]destination.west);

  \draw[line width=1.1pt, blue!70!black] (10.35,-1.25) -- (10.95,-1.25);
  \node[anchor=west, font=\scriptsize, text=blue!70!black]
    at (11.05,-1.25) {air coupling};
  \draw[line width=1.1pt, orange!85!black] (10.35,-1.65) -- (10.95,-1.65);
  \node[anchor=west, font=\scriptsize, text=orange!85!black]
    at (11.05,-1.65) {commuter coupling};
\end{tikzpicture}
\caption{Stochastic transfer on a directed multilayer network, specialized to
epidemic mobility. Population size $N_i$, gateway prominence $G_{i\ell}^{A}$, and
temporal fluctuation $X_{i\ell}^{A}(t)$ determine directed air and commuting
flows. Infectious prevalence then thins those flows into route-specific
infectious movements that can seed destination exposure.}
\label{fig:model_schematic}
\end{figure}

Figure~\ref{fig:model_schematic} shows the epidemic specialization used in the
numerical study. Air passengers may relocate between resident populations,
whereas commuters typically mix away from home and return. The reduced-form
experiment does not introduce explicit traveler or return compartments; it uses
layer-specific exposure couplings with different spatial kernels and transfer
scales. The theoretical results concern the stochastic flow-to-event map, while
the simulation embeds those events in a regional SEIR process.

\subsection{Directed spatial intensities}

For each $\ell$ and $t$, let $T_{\ell t}>0$ be the target conditional expected
layer volume. Let $W_{i\ell}^{A}(t)>0$ for
$A\in\{\mathrm{out},\mathrm{in}\}$, let
$\phi_\ell:[0,\infty)\to(0,\infty)$ be a measurable nonincreasing function
with $\phi_\ell(0)=1$, and let $r_{ij}^{(\ell)}\geq0$ be a time-invariant
dyadic affinity, with $r_{ii}^{(\ell)}=0$ and at least one positive off-diagonal
entry. Conditional on these quantities, take the off-diagonal counts to be
independent with
\begin{align}
A_{ij}^{(\ell)}(t)\mid\bm W(t)
&\sim \operatorname{Poisson}\{\mu_{ij}^{(\ell)}(t)\},\\
\mu_{ij}^{(\ell)}(t)
&=T_{\ell t}
\frac{W_{i\ell}^{\mathrm{out}}(t)W_{j\ell}^{\mathrm{in}}(t)
\phi_\ell(h_{ij})r_{ij}^{(\ell)}}
{\sum_{u\ne v}W_{u\ell}^{\mathrm{out}}(t)W_{v\ell}^{\mathrm{in}}(t)
\phi_\ell(h_{uv})r_{uv}^{(\ell)}},
\label{eq:mobility}
\end{align}
where $h_{ij}=\|\bm{x}_i-\bm{x}_j\|$. The denominator is assumed positive and
normalizes the intensities so that
$\sum_{i\ne j}\mu_{ij}^{(\ell)}(t)=T_{\ell t}$. Hence $T_{\ell t}$ is the
conditional mean total rather than the realized Poisson total; conditional on a
fixed realized total, the same construction gives a multinomial allocation over
ordered pairs. Monotonicity of $\phi_\ell$ encodes distance deterrence, while
$r_{ij}^{(\ell)}$ is the nonnegative multiplicative departure from the
node-weight--distance component for the ordered pair $(i,j)$. Only relative
values of the affinities are identified because their common scale cancels in
\eqref{eq:mobility}. Overdispersion can be represented by replacing the
conditional Poisson law with a negative-binomial observation model.

We factor node weight into common node scale, persistent residual structure,
and temporal fluctuation,
\begin{equation}
W_{i\ell}^{A}(t)
=\underbrace{N_i^{\beta_{\ell A}}G_{i\ell}^{A}}_{B_{i\ell}^{A}}
X_{i\ell}^{A}(t),
\qquad A\in\{\mathrm{out},\mathrm{in}\}.
\label{eq:decomp}
\end{equation}
Here $\beta_{\ell A}\in\mathbb R$ is the node-scaling exponent,
$G_{i\ell}^{A}>0$ is the time-invariant scale-residual factor, and
$X_{i\ell}^{A}(t)>0$ is the temporal multiplier. Thus
$B_{i\ell}^{A}:=N_i^{\beta_{\ell A}}G_{i\ell}^{A}$ is the persistent structural
weight. In the mobility application we call $G_{i\ell}^{A}$ \emph{gateway
prominence}. On the log scale, $\log G_{i\ell}^{A}$ is the residual after removing
$\beta_{\ell A}\log N_i$ from $\log B_{i\ell}^{A}$. A location convention such
as $\sum_i\log G_{i\ell}^{A}=0$ and
$\E\{\log X_{i\ell}^{A}(t)\}=0$ separates the two factors; any common
layer-direction scale is absorbed by the normalization in
\eqref{eq:mobility}. For a fixed direction $A$, cross-layer gateway alignment
is dependence between the paired vectors
$(G_{ia}^{A},G_{ic}^{A})$, $i\in V$. The full vector $\bm B_i$ therefore
includes population-mediated dependence, whereas $\bm G_i$ isolates dependence
between population-residual transport roles. The field $\bm X(t)$ contains the
remaining period-specific variation. The census commuting matrix informs
persistent structure, while weekly air traffic informs temporal fluctuation.

\subsection{Measures of interlayer dependence}

Fix a direction $A\in\{\mathrm{out},\mathrm{in}\}$ and suppress it from the
notation. Let $F_\ell$ be the marginal distribution function of
$D_{i\ell}$. The interlayer upper-tail dependence coefficient is
\begin{equation}
\lambda_U^{(a,c)}=\lim_{q\uparrow1}
\Prb\{D_{ic}>F_c^{-1}(q)\mid D_{ia}>F_a^{-1}(q)\}.
\label{eq:lambda}
\end{equation}
Here $F_\ell^{-1}(q)=\inf\{x:F_\ell(x)\ge q\}$, and the limit is taken when it
exists \citep{coles1999,coles2001}. Applied to raw strengths,
\eqref{eq:lambda} measures actual joint mobility concentration and includes
population-mediated dependence. Its population-conditioned analogue is estimated
by replacing $D_{i\ell}$ with the fitted residual factor $\widehat G_{i\ell}$
from \eqref{eq:decomp}; it measures alignment of gateway roles beyond
demographic scale. For finite $n$, we use two explicitly finite-sample
diagnostics. The rank association is
\begin{equation}
\rho_S=\operatorname{Corr}\{\widehat F_a(D_{ia}),
\widehat F_c(D_{ic})\},
\label{eq:rank_association}
\end{equation}
and the top-$q$ overlap count is
\begin{equation}
K_q=\sum_{i=1}^n
\mathbf 1\{\widehat F_a(D_{ia})>q,\ \widehat F_c(D_{ic})>q\}.
\label{eq:finite_overlap}
\end{equation}
Under independent uniformly permuted rankings and fixed marginal top-set sizes,
$K_q$ has the corresponding hypergeometric null distribution. At the available
$n=110$, these diagnostics are more stable than attempting to estimate the
limit in \eqref{eq:lambda}; we report them for both raw and
population-residual strengths.

\subsection{Stochastic event transfer and epidemic specialization}

Let $q_{i\ell}(t)>0$ be a possibly random layer-specific transfer multiplier.
Conditional on the directed flow and the multiplier, define propagated route
counts by
\begin{equation}
Z_{ij\ell}(t)\mid A_{ij}^{(\ell)}(t),q_{i\ell}(t)
\sim\operatorname{Poisson}\{q_{i\ell}(t)A_{ij}^{(\ell)}(t)\},
\label{eq:route_exports}
\end{equation}
conditionally independently over ordered routes. With
$Z_{i\ell}^{\mathrm{out}}(t)=\sum_{j\ne i}Z_{ij\ell}(t)$, Poisson
superposition gives the node-total counterpart
\begin{equation}
Z_{i\ell}^{\mathrm{out}}(t)\mid D_{i\ell}^{\mathrm{out}}(t),q_{i\ell}(t)
\sim\operatorname{Poisson}\{q_{i\ell}(t)D_{i\ell}^{\mathrm{out}}(t)\}.
\label{eq:exports}
\end{equation}
This is the stochastic network-to-event map analyzed in
Section~\ref{sec:theory}. In the epidemic
specialization, let $p_i(t)=I_i(t)/N_i$ be infectious prevalence and set
$q_{i\ell}(t)=c_\ell p_i(t)$, where $c_\ell$ contains the travel window,
ascertainment, and layer-specific transfer scale. This parameterization prevents
resident population from being counted once in total mobility and again in the
infectious pool. If
$M_{i\ell}=D_{i\ell}^{\mathrm{out}}/N_i$ is per-capita mobility, then the
node-total intensity becomes $c_\ell I_iM_{i\ell}$. A corresponding
discrete-time reduced-form SEIR transition makes the destination exposure map
explicit. Conditional on the state and mobility at time $t$, the expected new
exposures in region $i$ combine a local force of infection with a directed
import force from each layer,
\begin{equation}
M_i(t)=\beta_i\,S_i(t)\,p_i(t)
+\frac{S_i(t)}{N_i}\sum_{\ell\in\mathcal L}c_\ell
\sum_{j\ne i}A_{ji}^{(\ell)}(t)\,p_j(t),
\label{eq:foi}
\end{equation}
and the new exposures are a susceptible-capped Poisson leap
$Y_i(t)=\min\{S_i(t),\operatorname{Poisson}(M_i(t))\}$. The first term is local
transmission and each term in the sum is the layer-$\ell$ import force.
Equation~\eqref{eq:route_exports} isolates the random
infectious-movement input whose tail is studied in
Section~\ref{sec:theory}, and equation~\eqref{eq:foi}
shows how those inputs enter a full epidemic transition. Later epidemic size
also depends on progression, ancestry, susceptible depletion, and intervention.

\section{Aggregation principle}
\label{sec:theory}

The radial magnitude of a multivariate extreme is its total size, whereas its
angular composition records the relative contributions of layers or routes.
For a power-law tail, the index $\tau$ controls how quickly large values become
rare. The following results determine how these objects transform under three
operations: shared-factor generation of layer intensities, random diagonal
sampling with Poissonization, and nonnegative route aggregation. The first two
preserve the radial index under moment conditions, and the latter preserves that
index while generally making the original angular measure nonidentifiable.

In words, a distribution is multivariate regularly varying when its rescaled
large values concentrate on rays whose directions follow a fixed spectral
measure. That spectral measure is the layer and route composition of an extreme
event, and the transfer results track how each operation moves it. The
definitions that follow fix this idea precisely.

To formulate the transfer results, we first fix the standard notion of regular
variation.
\paragraph{Definition 1.}
A positive random variable $R$ is \emph{regularly varying with index}
$\tau>0$ if its survival function $\overline F_R(x)=\Prb(R>x)$ satisfies
\begin{equation}
\frac{\overline F_R(rx)}{\overline F_R(r)}\longrightarrow x^{-\tau}
\qquad\text{as }r\to\infty,\quad x>0.
\label{eq:univariate_rv_definition}
\end{equation}
Equivalently, there is a scaling function
$b\in\mathrm{RV}_{1/\tau}$ such that
$t\Prb\{R/b(t)>x\}\to x^{-\tau}$; see \citet[Chapter~3]{resnick2007}.

We next define multivariate regular variation (MRV) in the $M$-convergence sense, as
given in \citet{lindskog2014regular}.
\paragraph{Definition 2.}
For a nonnegative random vector $\bm Y\in[0,\infty)^d$, let
$\mathbb E_d=[0,\infty)^d\setminus\{\bm0\}$. Let $\mathbb M(\mathbb E_d)$ be
the class of Borel measures $\mu$ on $\mathbb E_d$ satisfying
$\mu\{\bm y:\|\bm y\|\ge r\}<\infty$ for every $r>0$, and let
\begin{equation}
\mathcal C_{\mathbb E_d}=\left\{f:\mathbb E_d\to[0,\infty):
\begin{array}{l}
f\text{ is bounded and continuous, and }f(\bm y)=0\\
\text{whenever }\|\bm y\|<r_f\text{ for some }r_f>0
\end{array}\right\}.
\label{eq:mrv_test_class}
\end{equation}
For $\mu_t,\mu\in\mathbb M(\mathbb E_d)$, write
$\mu_t\xrightarrow{M(\mathbb E_d)}\mu$ if
\begin{equation}
\int f\,d\mu_t\longrightarrow\int f\,d\mu
\qquad\text{for every }f\in\mathcal C_{\mathbb E_d}.
\label{eq:m_convergence_main}
\end{equation}
A distribution is MRV with index $\tau>0$ if, for some
$b\in\mathrm{RV}_{1/\tau}$ with $b(t)\to\infty$ and some nonzero
$\mu\in\mathbb M(\mathbb E_d)$,
\begin{equation}
t\Prb\{b(t)^{-1}\bm Y\in\cdot\}
\xrightarrow{M(\mathbb E_d)}\mu(\cdot),
\label{eq:mrv_definition_main}
\end{equation}
and the limit is homogeneous, $\mu(rA)=r^{-\tau}\mu(A)$ for $r>0$. After a
norm is fixed, its probability spectral measure is the normalized push-forward
of $\mu$ under $\bm y\mapsto\bm y/\|\bm y\|$ restricted to
$\{\bm y:\|\bm y\|>1\}$.

\subsection{Shared-factor generation of interlayer tail dependence}

The decomposition in \eqref{eq:decomp} suggests more than assuming regular
variation of the completed structural vector. A common heavy-tailed node factor
can generate that regular variation and yields an explicit interlayer tail
coefficient. Resident population is the common factor in the mobility
application, but the result is not tied to that interpretation. We use the
generalized Breiman theorem of
\citet[Appendix~A]{wang2024heterogeneous}. Unlike the classical
independent-product setting, the multiplier here may be evaluated at the
heavy-tailed radial variable itself: $\bm A(R)$ may change with node scale,
provided it converges at large scales and obeys a uniform moment bound.

The following result derives the stochastic-transfer limit from a common
heavy-tailed factor.
\begin{theorem}
\label{thm:shared_factor}
Let $R$ be positive and regularly varying with index $\tau>0$ and scaling
$b(t)$:
\[
t\Prb\{R/b(t)>x\}\longrightarrow x^{-\tau},\qquad x>0.
\]
Let $\{\bm A(r):r\geq0\}$ be an $m$-dimensional positive stochastic process,
independent of $R$, such that
\[
\bm A(r)\longrightarrow\bm A_\infty\quad\text{almost surely as }r\to\infty
\]
for a finite positive random vector $\bm A_\infty$. Suppose, for some
$\epsilon>0$ and some norm,
\begin{equation}
\sup_{r\geq0}\E\|\bm A(r)\|^{\tau+\epsilon}<\infty.
\label{eq:wr_uniform_moment}
\end{equation}
Conditional on $(R,\bm A)$, let $\bm Z$ have Poisson marginal means
\[
\bm\Lambda=R\bm A(R).
\]
Then $\bm\Lambda$ and $\bm Z$ are multivariate regularly varying with index
$\tau$, scaling $b(t)$, and the same limit measure
\begin{equation}
\mu_A(C)=
\E\left[
\int_0^\infty
\mathbf 1\{r\bm A_\infty\in C\}
\tau r^{-\tau-1}\,dr
\right].
\label{eq:shared_factor_measure}
\end{equation}
Set $m_k=\E(A_{\infty,k}^{\tau})$. For coordinates $j\ne k$, their
marginally standardized upper-tail dependence coefficient is
\begin{equation}
\lambda_{U}^{(j,k)}
=
\E\left[
\min\left\{
\frac{A_{\infty,j}^{\tau}}{m_j},
\frac{A_{\infty,k}^{\tau}}{m_k}
\right\}
\right].
\label{eq:shared_factor_lambda}
\end{equation}
\end{theorem}

Theorem~\ref{thm:shared_factor} is a \emph{generative} result: it starts from a
one-dimensional heavy-tailed primitive $R$ and proves regular variation of the
resulting vector, while also yielding the explicit coefficient
\eqref{eq:shared_factor_lambda}. Its process-indexed step follows
\citet[Appendix~A]{wang2024heterogeneous}; the complete specialization and the
Poissonization argument are given in the electronic supplementary material, section~S10.

For \eqref{eq:decomp} with a common node-scale exponent $\beta$, take
$R=N_i^\beta$. If $N_i$ has tail index $\alpha$, then
$R$ has index $\tau=\alpha/\beta$. The vector $\bm A(r)$ can contain
residual structure, transfer intensity, observation scale, and temporal
multipliers. A constant process recovers the usual assumption that these
factors are independent of node scale. The process formulation also allows
their distribution to vary with node scale and converge to a stable large-node
regime. Equation~\eqref{eq:shared_factor_lambda} separates the common radial
source from the limiting layer-specific composition. In the mobility
application, $N_i$ is population and $\bm A(r)$ contains gateway, prevalence,
and layer-transfer factors.
For unequal exponents $\beta_{\ell A}$, coordinate $k$ instead has marginal
tail index $\alpha/\beta_k$ under the corresponding moment condition; a
single standard radial scaling is then generally inappropriate, motivating
the marginal standardization used in the empirical comparisons.

\subsection{Structural transfer under random sampling}

The next result starts one stage later. Rather than deriving multivariate
regular variation from a common scalar factor, it assumes that the completed
structural vector $\bm B_i$ is already multivariate regularly varying and asks
how an independent, layer-specific sampling matrix transforms its limit
measure. It therefore covers structural dependence not representable by a
single common factor, but it does not by itself produce the explicit
tail-dependence formula in Theorem~\ref{thm:shared_factor}.

The next result transfers an assumed structural MRV limit through an independent
layer-specific sampling map.
\begin{theorem}
\label{thm:structural}
Fix node $i$ and let $\bm B_i=(B_{i\ell})_{\ell\in\mathcal L}$ be its full
structural layer vector, including common node scale, and assume it is
multivariate regularly varying on
$\mathbb E=[0,\infty)^{|\mathcal L|}\setminus\{\bm0\}$ with index $\tau>0$,
scaling function $b(t)$, and nonzero Radon limit measure $\nu$.
Let $Q_i$ be a random diagonal matrix, independent of $\bm B_i$, with positive
diagonal entries and
$\E\|Q_i\|^{\tau+\epsilon}<\infty$ for some $\epsilon>0$. If the entries of
$\bm Z_i$, conditional on $(\bm B_i,Q_i)$, have Poisson marginal distributions
with mean vector $Q_i\bm B_i$, then $\bm Z_i$ is multivariate regularly varying with
index $\tau$ and limit measure
\begin{equation}
\mu_Z(A)=\E\left[\nu\{\bm b:Q_i\bm b\in A\}\right]
\label{eq:random_pushforward}
\end{equation}
for continuity sets $A$ bounded away from the origin.
\end{theorem}

For the mobility specialization,
$Q_i=\operatorname{diag}\{c_\ell p_iX_{i\ell}\}_{\ell\in\mathcal L}$ combines
prevalence, layer transfer, and temporal fluctuation. The theorem itself applies
to any independent positive diagonal sampling map satisfying the stated moment
condition.

Note also that Theorems~\ref{thm:shared_factor} and \ref{thm:structural} start from different
objects. Theorem~\ref{thm:shared_factor} derives an MRV limit for
$R\bm A(R)$ from a common heavy-tailed factor $R$. Theorem~\ref{thm:structural}
instead assumes that the structural vector $\bm B_i$ is MRV and gives the limit
after the independent sampling map $Q_i$. In the common-factor setting, the
first theorem can supply the MRV input required by the second theorem.

Equation~\eqref{eq:random_pushforward} preserves the index $\tau$ because it
is a pushforward of a homogeneous limit measure. Its normalized spectral
measure, and hence $\lambda_U$, can nevertheless change when the diagonal
entries of $Q_i$ vary randomly by layer. A fixed diagonal $Q_i$ only rescales
coordinates, and a common scalar multiplier cancels after normalization.
Unequal exponents $\beta_{\ell A}$ can require coordinate-specific scaling, so
they need not yield one common radial index. Proof details are in
the electronic supplementary material, section~S10.

\subsection{Route resolution, aggregation, and source uncertainty}

Theorem~\ref{thm:structural} acts on the layer-total vector
$\bm B_i=(B_{i\ell})_{\ell\in\mathcal L}$ of one source node. To obtain
route-level quantities, collect the source-layer totals in a vector and apply
a deterministic allocation map before event sampling. This construction yields
the $|\mathcal R|$-dimensional route vector without a separate heavy-tail
assumption for every route.
Formally, let $\mathcal A$ index a finite set of source-layer cells
$a=(i,\ell)$, and let
$\bm B^{A}=(B_a:a\in\mathcal A)$ contain their structural movement totals.
Let $\mathcal R$ be a finite set of directed route-layer coordinates
$e=(i,j,\ell)$, and let $P$ be a fixed nonnegative
$|\mathcal R|\times|\mathcal A|$ route-allocation matrix. Its columns sum to
one, and $P_{e,a}=0$ unless route $e$ leaves the source-layer cell $a$.
Thus $P\bm B^{A}$ distributes each source-layer total over destinations. Write
$\bm Z^{R}=(Z_e:e\in\mathcal R)$ for propagated route counts. This
representation retains source, destination, and layer identity; source totals,
destination totals, and layer totals are nonnegative
linear aggregations of $\bm Z^{R}$. In \eqref{eq:mobility}, a column of $P$
is the corresponding set of destination weights normalized within a fixed
source and layer.

The following proposition extends the transfer map from source-layer totals to
route-level counts and their dominant-route probabilities.
\begin{proposition}
\label{prop:routes}
Suppose $\bm B^{A}$ is multivariate regularly varying with index $\tau$,
scaling $b(t)$, and limit measure $\nu_A$. Let $Q^{R}$ be a random positive
diagonal matrix, independent of $\bm B^{A}$, such that
$\E\|Q^{R}\|^{\tau+\epsilon}<\infty$ for some $\epsilon>0$.
Conditional on $(\bm B^{A},Q^{R})$, suppose each $Z_e$ has Poisson mean
$(Q^{R}P\bm B^{A})_e$. Then $\bm Z^{R}$ is multivariate regularly varying
with index $\tau$ and limit measure
\begin{equation}
\mu_R(A)=\E\left[\nu_A\{b:Q^{R}Pb\in A\}\right].
\label{eq:route_limit}
\end{equation}
\begin{enumerate}
\item Let $S$ be a fixed nonnegative matrix with no zero column. Then
$S\bm Z^{R}$ is multivariate regularly varying with index $\tau$ and limit
measure. This aggregation retains the radial index but can map distinct
route-level angular measures to the same aggregate measure:
\begin{equation}
\mu_S(A)=\mu_R\{z:Sz\in A\},
\label{eq:route_pushforward}
\end{equation}
\item Use the $\ell_1$ norm, write
$\mathbb S_+=\{w\geq0:\|w\|_1=1\}$, and let $H_R$ be the probability spectral
measure obtained by normalizing $\mu_R$ on $\{z:\|z\|_1>1\}$. For
\[
C_e=\{w\in[0,1]^{|\mathcal R|}:\|w\|_1=1,
\ w_e>\max_{f\ne e}w_f\},
\]
assume that $H_R$ assigns zero mass to directions tied for the largest
coordinate. Let $J(z)$ return the coordinate of the largest entry of $z$,
using any fixed deterministic rule when the maximum is tied. Then
\begin{equation}
\Prb\!\left\{J(\bm Z^{R})=e
\,\middle|\, \|\bm Z^{R}\|_1>u\right\}
\longrightarrow H_R(C_e),
\label{eq:dominant_route}
\end{equation}
Consequently, the entropy of the conditional dominant-route distribution
converges to
\begin{equation}
\mathcal H_R=-\sum_{e\in\mathcal R}H_R(C_e)\log H_R(C_e),
\label{eq:route_entropy}
\end{equation}
with $0\log0=0$.
\end{enumerate}
\end{proposition}

Equation~\eqref{eq:dominant_route} identifies the route most likely to carry
the largest realized propagated count when total activity is extreme. Its
angular mass records the source-destination-layer coordinate that dominates
an unusually active system. In the epidemic application this is a
tail-specific counterpart to route-occurrence and source-entropy summaries in
stochastic reconstructions of early spread \citep{zhang2026spatial}; unlike
unconditional ensemble occurrence, it conditions on extreme total infectious
movement.

The following corollary gives the spectral measure induced by a nonnegative
aggregation map and states the resulting nonidentifiability condition.
\begin{corollary}
\label{cor:aggregation_spectral}
Under Proposition~\ref{prop:routes}, let $H_R$ be the route-level probability
spectral measure under the $\ell_1$ norm. For a nonnegative aggregation
matrix $S$ with no zero column, the probability spectral measure of
$S\bm Z^R$ is
\begin{equation}
H_S(A)=
\frac{
\displaystyle\int_{\mathbb S_+}
\|Sw\|_1^\tau
\mathbf 1\left\{\frac{Sw}{\|Sw\|_1}\in A\right\}H_R(dw)}
{\displaystyle\int_{\mathbb S_+}\|Sw\|_1^\tau H_R(dw)}.
\label{eq:aggregate_spectral}
\end{equation}
If the normalized map
$T_S(w)=Sw/\|Sw\|_1$ is not injective on $\mathbb S_+$, then there exist two
distinct route-level spectral measures
$H_R^{(1)}\ne H_R^{(2)}$ that induce the same aggregate spectral measure
$H_S$.
\end{corollary}

The weight $\|Sw\|_1^\tau$ in \eqref{eq:aggregate_spectral} scales each route
composition by its retained aggregate magnitude before normalization.
Nonidentifiability occurs whenever two route compositions lie on the same
normalized aggregation fiber. The two
allocation matrices used in the finite-sample experiment give a construction
of this phenomenon; implementation details and the corresponding figure are
provided in the electronic supplementary material, section~S2.


The next proposition characterizes when aggregation preserves enough
information to recover route composition.
\begin{proposition}
\label{prop:identifiability}
Let $S$ be a nonnegative $m\times|\mathcal R|$ matrix with no zero column, let
$s_e$ be its $e$th column, and define
$\hat s_e=s_e/\|s_e\|_1$. Thus $\hat s_e$ is the normalized vector of retained
aggregate coordinates generated by a unit contribution on route $e$. The map
$T_S(w)=Sw/\|Sw\|_1$ is injective on $\mathbb S_+$ if and only if the signatures
$\{\hat s_e:e\in\mathcal R\}$ are affinely independent. When this holds:
\begin{enumerate}
\item (Counting.) Injectivity requires $|\mathcal R|\le m$. Full pooling to a
single aggregate coordinate ($m=1$) is injective only if $|\mathcal R|=1$.
\item (Recovery.) The route-level spectral measure is identified from the
aggregate spectral measure: writing $\psi=T_S^{-1}$, for every bounded
continuous $\varphi$ on $\mathbb S_+$,
\begin{equation}
\int\varphi\,dH_R=
\frac{\displaystyle\int
\varphi(\psi(v))\,\|S\psi(v)\|_1^{-\tau}\,H_S(dv)}
{\displaystyle\int \|S\psi(v)\|_1^{-\tau}\,H_S(dv)}.
\label{eq:spectral_recovery}
\end{equation}
\end{enumerate}
For layer pooling, let $S$ sum route coordinates within each
source-destination pair and write
$\bar w_{ij}=(Sw)_{ij}=\sum_\ell w_{ij\ell}$. Then $\|Sw\|_1=1$ on
$\mathbb S_+$, and the aggregate spectral measure is exactly the
source-destination pushforward
\begin{equation}
H_{\rm OD}(A)=H_R\{w:Sw\in A\}.
\label{eq:od_pushforward}
\end{equation}
Thus layer pooling identifies the source-destination spectral measure. Routes
that share $(i,j)$ and differ only in $\ell$ have identical signatures, so their
within-pair layer composition remains unidentified.
\end{proposition}

The proposition gives a geometric identifiability criterion. A route
composition $\bm w$ is first mapped to $S\bm w$ and then normalized. This map is
one-to-one exactly when the route signatures are affinely independent, which
requires at least $|\mathcal R|$ retained aggregate coordinates. A destination's
total importation remains observable because it is an aggregate coordinate.
Mode attribution is lost when layer-specific routes share one
source-destination signature. Full pooling ($m=1$) therefore loses all route
composition whenever more than one route is present. The proof is in
the electronic supplementary material, section~S10.

The following result connects an extreme route-intensity vector to the route
and time of the first established importation.
\begin{theorem}
\label{thm:first_invasion}
Let $\bm\Lambda=(\Lambda_e:e\in\mathcal R)$ be a positive route-intensity
vector that is multivariate regularly varying with index $\tau$, scaling
$b(t)$, and probability spectral measure $H_R$ under the $\ell_1$ norm.
For route $e$, let $s_e\in[0,1]$ be the probability that a propagated event on
that route is successful. Conditional on $\bm\Lambda$, suppose successful
events on the routes are independent homogeneous Poisson processes with rates
$s_e\Lambda_e$. Let $T_*$ be the time of the first successful event and $E_*$
its route. Write
$\bm s=(s_e:e\in\mathcal R)$ and let $\bm w\in\mathbb S_+$ denote
a spectral direction, and define
$\bm s^\top\bm w=\sum_{e\in\mathcal R}s_ew_e$. Assume
\begin{equation}
H_R\{\bm w:\bm s^\top\bm w=0\}=0.
\label{eq:positive_establishment_mass}
\end{equation}
Then, for every route $e$,
\begin{equation}
\Prb\{E_*=e\mid\|\bm\Lambda\|_1>b(t)\}
\longrightarrow
\pi_e^*
:=\int_{\mathbb S_+}
\frac{s_ew_e}{\bm s^\top\bm w}\,H_R(d\bm w).
\label{eq:first_route_limit}
\end{equation}
More generally, for every $x\geq0$,
\begin{align}
&\Prb\{E_*=e,\ b(t)T_*>x
\mid\|\bm\Lambda\|_1>b(t)\}\nonumber\\
&\quad\longrightarrow
\int_{\mathbb S_+}\int_1^\infty
\frac{s_ew_e}{\bm s^\top\bm w}
\exp\{-xr\,\bm s^\top\bm w\}
\tau r^{-\tau-1}\,dr\,H_R(d\bm w),
\label{eq:first_route_time_limit}\\
&\Prb\{b(t)T_*>x\mid\|\bm\Lambda\|_1>b(t)\}\nonumber\\
&\quad\longrightarrow
\int_{\mathbb S_+}\int_1^\infty
\exp\{-xr\,\bm s^\top\bm w\}
\tau r^{-\tau-1}\,dr\,H_R(d\bm w).
\label{eq:first_time_limit}
\end{align}
Consequently, the first-success route entropy converges to
\begin{equation}
\mathcal H_*=-\sum_{e\in\mathcal R}\pi_e^*\log\pi_e^*.
\label{eq:first_invasion_entropy}
\end{equation}
\end{theorem}

The conditional Poisson-process statement specifies competing clocks given
$\bm\Lambda$; after integrating over $\bm\Lambda$, the limits in
\eqref{eq:first_route_limit}-\eqref{eq:first_time_limit} are conditioned only
on the extreme-intensity event. There is therefore no second, independent tail
conditioning. The success probabilities can represent route reliability,
detection, activation, or survival of an initiated event. In the epidemic
specialization, a local branching approximation gives the establishment
probability: if route $e$ ends at destination $j(e)$ with offspring generating
function $G_{j(e)}$, then $s_e=1-q_{j(e)}$, where $q_{j(e)}$ is the smallest
solution of $G_{j(e)}(q)=q$.

Theorem~\ref{thm:first_invasion} differs from
\eqref{eq:dominant_route} in both the event being selected and its weighting.
Equation~\eqref{eq:dominant_route} selects the route carrying the largest
realized propagated count under an extreme total. Theorem~\ref{thm:first_invasion}
instead selects the route on which the first successful event occurs. Its
extremal share $w_e$ is reweighted by $s_e$ and normalized by
$\bm s^\top\bm w$; a lower-intensity route can win when its events are more
likely to succeed. Equations
\eqref{eq:first_route_time_limit}-\eqref{eq:first_time_limit} connect the
radial magnitude of an extreme to the accelerated first-event time scale. For
epidemic invasion, ``success'' means that an imported infection establishes a
self-sustaining local chain. A direct simulation of a model that is only
asymptotically regularly varying confirms all three predictions: as the
conditioning threshold grows, the empirical first-success route distribution,
its entropy, and the rescaled first-time law converge to
\eqref{eq:first_route_limit}, \eqref{eq:first_invasion_entropy}, and
\eqref{eq:first_time_limit}, and the first-success distribution is numerically
distinct from the dominant-route distribution
(the electronic supplementary material, section~S8).

The next corollary extends the transfer results to independent Gamma mixing and
negative-binomial event sampling.
\begin{corollary}
\label{cor:nb}
Let $(\bm B_*,M,\nu_*)$ denote a generic input vector, sampling map, and MRV
limit measure. For Theorem~\ref{thm:structural}, take
$(\bm B_*,M,\nu_*)=(\bm B,Q,\nu)$. For Proposition~\ref{prop:routes}, take
$(\bm B_*,M,\nu_*)=(\bm B^{A},Q^{R}P,\nu_A)$.
In either case, replace the
Poisson sampling mechanism by
\[
Z_k\mid \bm B_*,M,G\sim
\operatorname{Poisson}\!\left\{G_k(M\bm B_*)_k\right\},
\]
where $M$ has $m$ rows and
$G=\operatorname{diag}(G_1,\ldots,G_m)$ is positive and independent
of $(\bm B_*,M)$. Here $\bm B_*$ is the input vector, $M$ maps it into the
$m$ propagated-event coordinates, and $G_k$ is the sampling multiplier for
coordinate $k$. If
$\E\|GM\|^{\tau+\delta}<\infty$ for some $\delta>0$, the
conclusions hold with limit measure
\begin{equation}
\mu_{GM}(A)=\E\left[\nu_*\{b:GMb\in A\}\right].
\label{eq:nb_pushforward}
\end{equation}
In particular, $G_k\sim\operatorname{Gamma}(\kappa_k,\kappa_k)$
gives a negative-binomial marginal with conditional mean $(M\bm B_*)_k$
and variance
$(M\bm B_*)_k+(M\bm B_*)_k^2/\kappa_k$. If $G_k=G_0$ for every $k$,
the base limit measure is multiplied by $\E(G_0^\tau)$ and its normalized
spectral measure is unchanged;
coordinate-specific Gamma multipliers generally change that spectral measure.
\end{corollary}

Since Gamma variables have moments of every positive order, their use adds
no moment restriction beyond the corresponding assumption on $M$.
The corollary incorporates Gamma mixing without changing the structural tail
index. The Gamma factor can represent reporting, transfer, or event-generation
heterogeneity; in the epidemic specialization it can encode effective seeding
or first-generation offspring. A common multiplier changes radial magnitude,
whereas route- or layer-specific heterogeneity can also change the composition
of an extreme event.
Proofs of Proposition~\ref{prop:routes} and Corollary~\ref{cor:nb} are in
the electronic supplementary material, section~S10.

\section{Numerical and empirical evidence}
\label{sec:evidence}

\subsection{Synthetic epidemic study}
\label{sec:simulation}


The simulation study links the asymptotic transfer results to finite networks
and traces their consequences in a stochastic metapopulation epidemic. The
design is selective disruption: each mechanism-removal comparison alters exactly
one dependence channel, holds the other components and the total layer volume
fixed, and reads off the paired difference. Two constructions then test the
aggregation claim directly. An observational-equivalence twin builds two
multiplexes that share one pooled network, produce the identical epidemic, and
split that network into air and commuting differently; a population-tail
sensitivity sweep regenerates the network across four values of the population
tail parameter. All results use fixed seeds and are reproduced by the supplied
scripts. Finite-sample checks of
radial-index preservation, gateway residualization, tail-dependence estimation,
and route-level nonidentifiability are reported in the electronic supplementary
material, sections~S1--S2.

\subsubsection{Synthetic multiplex and epidemic model}
\label{sec:sim_model}

\paragraph{Directed multiplex generator.}

Each network has $n=60$ regions at independent uniform coordinates
$u_i\in[0,1]^2$, with Euclidean separation $d_{ij}=\lVert u_i-u_j\rVert$. Region
$i$ has resident population
\begin{equation}
N_i=40{,}000\,(1+P_i),\qquad P_i\sim\mathrm{Pareto}(\alpha),\quad\alpha=1.8,
\label{eq:sim_population}
\end{equation}
truncated at $900{,}000$; write $n_i=N_i/\operatorname{median}_kN_k$. For each
direction $A\in\{\mathrm{out},\mathrm{in}\}$ we draw standard normals with
cross-layer correlation $0.45$ and set the log-gateway factors
$\log G^{A}_{i\ell}=0.65\,z^{A}_{i\ell}$, so each gateway is lognormal with
log-standard-deviation $0.65$. Directed emission and attraction weights combine
demographic scale with the gateway factor,
\begin{equation}
w^{A}_{i,\mathrm{air}}=n_i^{0.95}\,G^{A}_{i,\mathrm{air}},\qquad
w^{A}_{i,\mathrm{commute}}=n_i^{1.05}\,G^{A}_{i,\mathrm{commute}},
\end{equation}
and the mean directed flow is a gravity form with zero diagonal,
\begin{equation}
\bar A^{\ell}_{ij}\propto w^{\mathrm{out}}_{i\ell}\,w^{\mathrm{in}}_{j\ell}\,
e^{-d_{ij}/\rho_\ell},\qquad
\rho_{\mathrm{air}}=0.65,\ \rho_{\mathrm{commute}}=0.13,
\end{equation}
rescaled so that $\sum_{ij}\bar A^{\mathrm{air}}_{ij}=0.002\sum_iN_i$ and
$\sum_{ij}\bar A^{\mathrm{commute}}_{ij}=0.075\sum_iN_i$ movements per day. The
slower air kernel decay ($\rho_{\mathrm{air}}>\rho_{\mathrm{commute}}$) gives
aviation the longer reach.

A daily unit-mean lognormal field modulates each layer,
$x_{i\ell}(t)=\exp\{\sigma\eta_{i\ell}(t)-\tfrac12\sigma^2\}$ with $\sigma=0.38$
and
\begin{equation}
\eta_{i\ell}(t)=\sqrt{f_\ell}\,c_\ell(t)+\sqrt{1-f_\ell}\,s_{i\ell}(t).
\end{equation}
Here $c_\ell(t)$ is a spatially constant common mode, $s_{i\ell}(t)$ is a
mean-zero Gaussian spatial field with correlation $e^{-d_{ij}/r_\ell}$, and the
common-mode fraction is $f_{\mathrm{air}}=0.40$, $f_{\mathrm{commute}}=0.25$.
The spatial ranges are $r_{\mathrm{air}}=0.35$, $r_{\mathrm{commute}}=0.18$, and
both the common modes and the spatial innovations carry cross-layer correlation
$0.35$. The realized day-$t$ flow scales each edge by the geometric mean of its
endpoint fields, $A^{\ell}_{ij}(t)=\bar A^{\ell}_{ij}\sqrt{x_{i\ell}(t)\,x_{j\ell}(t)}$.

\paragraph{Stochastic SEIR dynamics and route attribution.}

Each region carries susceptible, exposed, infectious, and removed counts
$S_i,E_i,I_i,R_i$ summing to $N_i$, with infectious prevalence
$p_i(t)=I_i(t)/N_i$. The epidemic runs for $110$ days from a seed of $20$
infectious individuals at the largest-outflow air source. New exposures follow
the exposure map \eqref{eq:foi} with a constant local coefficient
$\beta_i\equiv\beta=0.31$ and layer transfer scales
$\kappa_{\mathrm{air}}=0.075$, $\kappa_{\mathrm{commute}}=0.012$ in place of
$c_\ell$: region $j$ gains
$Y_j(t)=\min\{S_j(t),\operatorname{Poisson}(M_j(t))\}$ exposures on day $t$.
Progression is geometric: the daily
$E\to I$ and $I\to R$ transitions are $\operatorname{Binomial}(E_j,1-e^{-1/3})$
and $\operatorname{Binomial}(I_j,1-e^{-1/5})$, giving mean latent and infectious
periods of three and five days. Region $j$ is invaded on the first day that
$E_j+I_j+R_j\ge5$, and secondary invasion is the first such event outside the
seed. These simulation definitions differ from the reported-incidence onset used
in the real-data analysis.

Each daily exposure count is split by Poisson superposition, first among its
local, air, and commuting components in proportion to their means in
\eqref{eq:foi}, then among origins in proportion to the route means
$\kappa_\ell A^{\ell}_{ij}(t)p_i(t)$. A separate random stream carries this
attribution, so recording routes leaves the paired epidemic paths unchanged. The
route ledger of a replicate records counts on each directed
source-destination-layer route $\rho$; let
$\bm w=(\pi_\rho)$, with $\pi_\rho$ the route's share of all imported exposures,
be the ledger's angular component $\bm Z^{R}/\lVert\bm Z^{R}\rVert_1$. Conditional
on an extreme outbreak, $\bm w$ has the route-level spectral law $H_R$ of
Proposition~\ref{prop:routes}, and we summarize route concentration by two
bounded-continuous functionals of $H_R$: the route entropy
$\int(-\sum_\rho w_\rho\log w_\rho)\,H_R(d\bm w)$ (in nats) and the
dominant-route share $\int\max_\rho w_\rho\,H_R(d\bm w)$. We estimate both by
conditioning on the replicates whose total imported exposures fall in the upper
quartile, each system on its own outbreak-size tail, and averaging the
angular entropy and the top-route share over those replicates. This angular
functional is the informative spectral summary here: because the seed is fixed
at one hub, a single route is the largest coordinate in essentially every
replicate, so the entropy of the dominant-route ($\arg\max$) distribution
in \eqref{eq:route_entropy} is degenerate at the epidemic level and is validated
instead on the finite-sample route vector in the electronic supplementary material, section~S8.

\paragraph{Mechanism-removal comparisons.}

Each comparison modifies one channel and rescales the affected layer to its
original total. Removing population alignment permutes the demographic scale
$n_i$ in the commuting weights and leaves air fixed; removing gateway alignment
permutes the commuting gateway factors $G_{\mathrm{commute}}$. Randomizing
spatial placement applies one daily node permutation to both temporal fields,
stripping geographic placement while keeping the cross-layer pairing; the
spatial-plus-interlayer variant uses independent daily permutations per layer.
Two aggregation ablations, used in Section~\ref{sec:aggregation_hides}, complete
the set. The first sets $\kappa_{\mathrm{air}}$ and
$\kappa_{\mathrm{commute}}$ to their common volume-weighted value. The
composite static aggregate additionally sums the two mean flows and removes the
temporal fields. Three hundred paired replicates use common random numbers, and
paired bootstraps give the reported intervals.

\subsubsection{Consequences for epidemic invasion}

\paragraph{Invasion timing and spread.}

\begin{figure}[h]
\centering
\includegraphics[width=\textwidth]{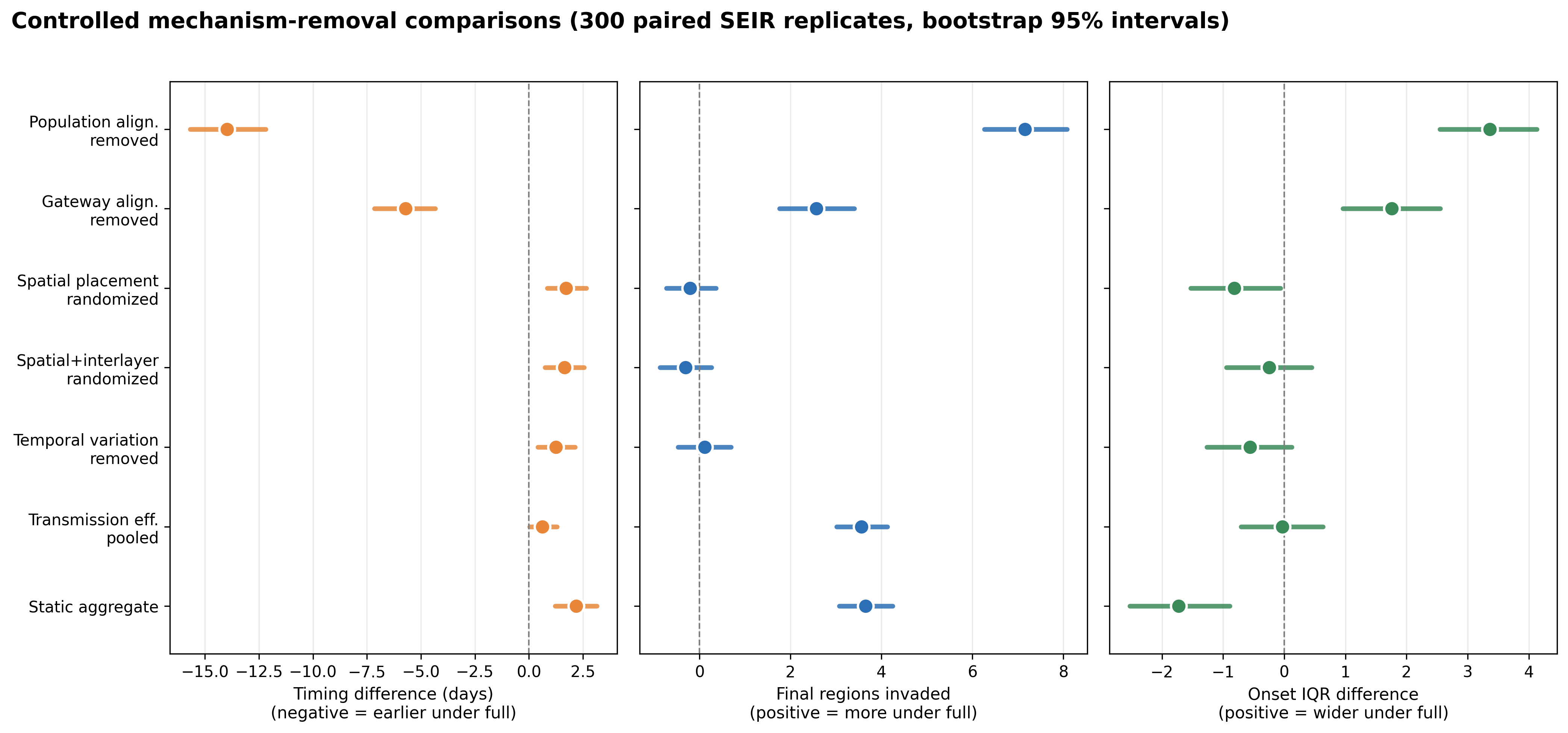}
\caption{Controlled mechanism-removal comparisons across 300 paired stochastic
SEIR replicates. Points are mean differences between the full multiplex and each
modified network, with paired bootstrap 95\% intervals. In the left panel a
negative value means secondary invasion occurs earlier under the full multiplex,
in the center a positive value means more regions are invaded, and the right
panel reports onset-time dispersion, where a positive value denotes a wider
interquartile range.}
\label{fig:epidemic_sim}
\end{figure}

Figure~\ref{fig:epidemic_sim} reports the mechanism-removal results for invasion
timing, spread, and onset-time dispersion.
In the baseline network the full multiplex reaches median first secondary
invasion on day 35, four regions by day 60, and 35 regions by day 110. Removing
population alignment delays secondary invasion by $13.98$ days (bootstrap 95\%
interval $12.14$--$15.85$), reduces day-60 spread by $1.72$ regions
($1.47$--$1.96$), and reduces day-110 spread by $7.16$ regions. Removing gateway
alignment delays invasion by $5.71$ days ($4.32$--$7.18$) and reduces day-60 and
day-110 spread by $0.49$ and $2.56$ regions.

The timing acceleration depends on a hub seed in this calibration. Across 12
independently generated networks, aviation, commuting, and population hub seeds
reproduce the effect; median-population and random seeds give near-zero timing
differences and slightly smaller final sizes (the electronic supplementary material, section~S3).
Route composition, which pooling governs, does not depend on timing. Temporal
spatial ablations are small in this calibration, and the electronic
supplementary material, section~S4 reports their common-mode sensitivity.

\paragraph{The infectious-route backbone.}

\begin{figure}[h]
\centering
\includegraphics[width=0.66\textwidth]{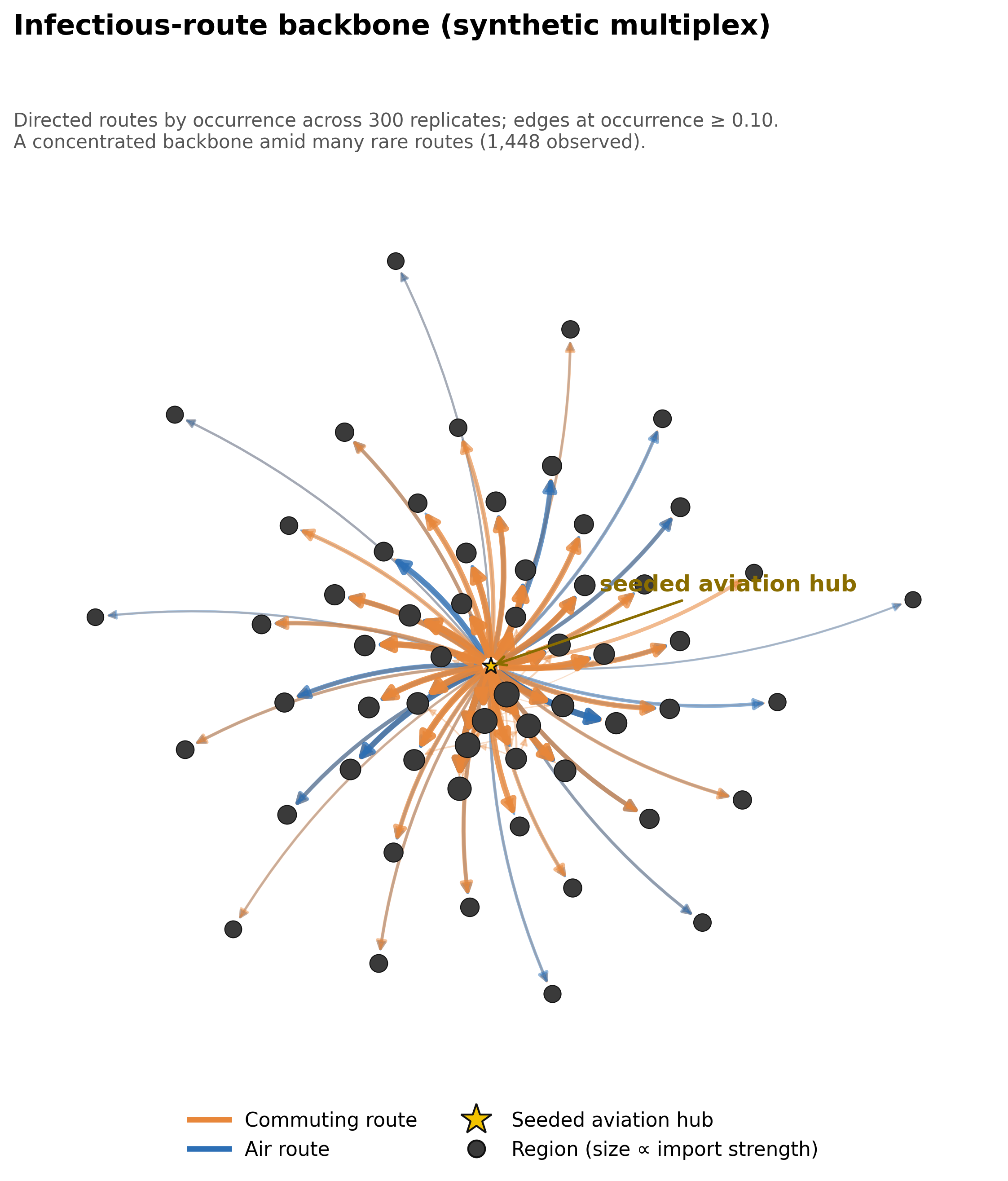}
\caption{Infectious-route backbone of the synthetic multiplex, shown as a
directed network over the 60 regions. Edges are infectious routes with ensemble
occurrence at least $0.10$ across 300 full-multiplex replicates, coloured orange
for commuting and blue for air, with width and opacity proportional to
occurrence probability. Node area is proportional to import strength, and the
seeded aviation hub is marked by a gold star.}
\label{fig:epidemic_routes}
\end{figure}

Figure~\ref{fig:epidemic_routes} displays the infectious-route backbone across
the 60 regions. Route attribution reveals substantial
pathway uncertainty beneath the aggregate outcomes. In the full multiplex, the
median replicate contains $87.5$ active infectious routes and 454 imported
exposures. Conditioned on the upper-quartile outbreaks, the route entropy is
$3.06$ (bootstrap 95\% interval $3.03$--$3.10$) and the dominant-route share is
$0.395$ ($0.390$--$0.401$). Across replicates, $91.9\%$ of the 1,448 routes ever
observed occur in fewer than $20\%$ of simulations; $2.4\%$ occur in at least
$80\%$. The epidemic therefore combines a concentrated core of frequent routes
with many unstable peripheral routes, paralleling the ensemble uncertainty found
in reconstructed pandemic networks \citep{zhang2026spatial}.

Demographic alignment concentrates the infectious-route backbone. Relative to
the population-alignment comparison, the full multiplex raises the dominant-route
share by $0.247$ (paired bootstrap 95\% interval $0.228$--$0.265$) and lowers
route entropy by $0.746$ ($0.678$--$0.811$). Gateway alignment, by contrast, does
not detectably change the extreme-outbreak route composition: its share and
entropy differences are $0.024$ ($-0.008$--$0.057$) and $-0.007$
($-0.114$--$0.100$), both intervals covering zero and both varying in sign across
conditioning thresholds. The gateway channel therefore shifts invasion timing
without concentrating the tail-conditioned route backbone. The two temporal
spatial comparisons change route entropy by at most $0.03$. Faster invasion in
the population-aligned multiplex is therefore associated with a more concentrated
principal route, and most realized routes remain rare across epidemic replicates.

\paragraph{Population-tail sensitivity.}

\begin{figure}[h]
\centering
\includegraphics[width=\textwidth]{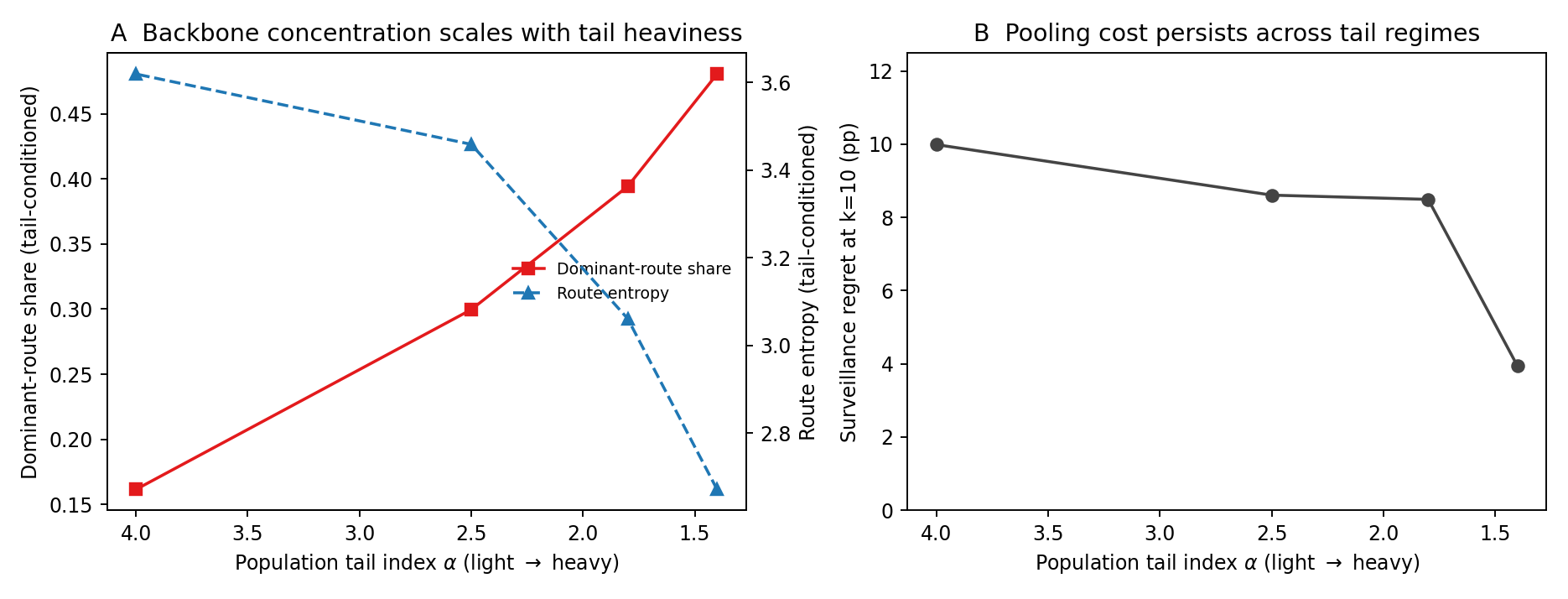}
\caption{Population-tail sensitivity of the synthetic epidemic, over 120
full-multiplex replicates at each value of $\alpha$ (greater generated
heterogeneity to the right). (A) Tail-conditioned dominant-route share and route
entropy against $\alpha$. (B) Pooled-versus-layer-resolved surveillance regret
$\mathcal R_{10}$ of \eqref{eq:regret} against $\alpha$.}
\label{fig:tail_sweep}
\end{figure}

Figure~\ref{fig:tail_sweep} reports population-tail sensitivity; Supplemental
Material Section~S9 gives design details and full results. In panel A, as the
population-tail parameter $\alpha$ falls from $4.0$
to $1.4$, increasing population heterogeneity in this generator, the
tail-conditioned dominant-route share rises from $0.162$ to $0.480$ and the
tail-conditioned route entropy falls from $3.62$ to $2.67$. Changing $\alpha$
regenerates the network and therefore changes hub dominance, route allocation,
and epidemic dynamics together. The trend shows that greater population
heterogeneity concentrates imported exposure on fewer routes in this calibration.
The sweep does not isolate the radial index and does not validate the asymptotic
radial-transfer result.

Figure~\ref{fig:tail_sweep}B shows that the
surveillance cost of pooling, quantified in
Section~\ref{sec:aggregation_hides}, stays between $4$ and $10$ percentage points
across these settings and is largest at $\alpha=4.0$. When heterogeneity is
greatest, a few dominant hubs lead both layers, so the pooled ranking more
closely tracks the air ranking. The decision cost therefore persists across the
four population-tail settings.

\subsubsection{Aggregation hides the import pathway}
\label{sec:aggregation_hides}

\paragraph{Two aggregation operations and direct pooling.}

We separate two operations that are often both called aggregation.
\emph{Direct pooling} sums the layer-specific import counts while holding
the realized epidemic path fixed; it preserves each destination's total imports
by construction but removes the layer label on every import. The
\emph{composite static aggregate} goes further, summing the two mean flow
matrices and applying a single volume-weighted transmission coefficient, so it
removes temporal variation, layer identity, and layer-specific transmission
efficiencies at once. Decomposing that composite gap is clean: the full model's
tail-conditioned route entropy exceeds the static aggregate by $0.402$. Setting
the transmission coefficients to their common volume-weighted value accounts
for $0.393$, with temporal
averaging and label pooling making up the small remainder. The composite aggregate is
useful as a sensitivity analysis but does not isolate the effect of pooling
labels; direct pooling does.

\begin{table}[h]
\centering\small
\caption{Aggregation analysis in the baseline synthetic epidemic, contrasting
direct pooling with the composite static aggregate as defined in the text.}
\label{tab:aggregation_audit}
\begin{tabular}{p{0.71\textwidth}r}
\toprule
Quantity & Value\\
\midrule
Destinations whose total imports are retained by direct pooling & 60 of 60\\
Destinations with air-dominated imports before pooling & 14 of 60 (23\%)\\
Network-wide air share of infectious imports & 14\%\\
Top-five import destinations shared by the composite static aggregate & 5 of 5\\
Top-ten import destinations shared by the composite static aggregate & 9 of 10 (Jaccard $0.82$)\\
\bottomrule
\end{tabular}
\end{table}

Direct pooling answers the surveillance question. Which destinations receive
many imports, and through which layer do they receive them? The first target is
available from total pooled imports. The second is not. In the layer-resolved
simulation, air is the dominant infectious-import layer for 14 of 60
destinations ($23\%$), even though it supplies only $14\%$ of infectious imports
network-wide, as Table~\ref{tab:aggregation_audit} reports. Pooling preserves the
total import count for those destinations, but it cannot determine whether air or
commuting supplied that count. That discarded layer label is exactly the
information a mode-specific intervention would require.

\paragraph{The decision cost of pooling.}

The size of that loss is a decision cost within the simulation. Consider a
planner who can place a mode-specific air surveillance resource at $k$
destinations. Let $a_j$ be destination $j$'s air-import burden, its total
air-layer infectious imports, and let $t_j=a_j+c_j$ be its total imports across
the air and commuting layers. Write $\mathcal O_k$ for the set of $k$
destinations with the largest $a_j$ and $\mathcal P_k$ for the set of $k$ with
the largest $t_j$, and let
$C(\mathcal S)=\sum_{j\in\mathcal S}a_j\big/\sum_j a_j$ be the fraction of
network air-import burden captured by a set $\mathcal S$. With layer-resolved
counts the planner selects $\mathcal O_k$; with only pooled mobility the air
layer is unobserved, so the best available proxy is $\mathcal P_k$. The regret at
budget $k$ is the air coverage forgone by that substitution,
\begin{equation}
\mathcal R_k=C(\mathcal O_k)-C(\mathcal P_k).
\label{eq:regret}
\end{equation}
At $k=10$, $C(\mathcal O_{10})=46.5\%$ and $C(\mathcal P_{10})=38.3\%$, so
$\mathcal R_{10}=8.2$ percentage points (paired bootstrap 95\% interval
$7.6$--$8.9$); the pooled decision thus forfeits about one-fifth of the air
coverage attainable at the same budget. The misdirection is systematic: none of
the ten destinations a pooled planner would select are air-dominated, because
the largest total-import destinations are commuting hubs, so the pooled resource
is placed where air is not the operative pathway. The regret is $\mathcal R_k=6.8$,
$8.2$, $7.0$, and $9.2$ percentage points at $k=5,10,15,20$ (Supplemental
Material Section~S7).

\paragraph{Nonidentifiability: an observational-equivalence twin.}

\begin{figure}[h]
\centering
\includegraphics[width=\textwidth]{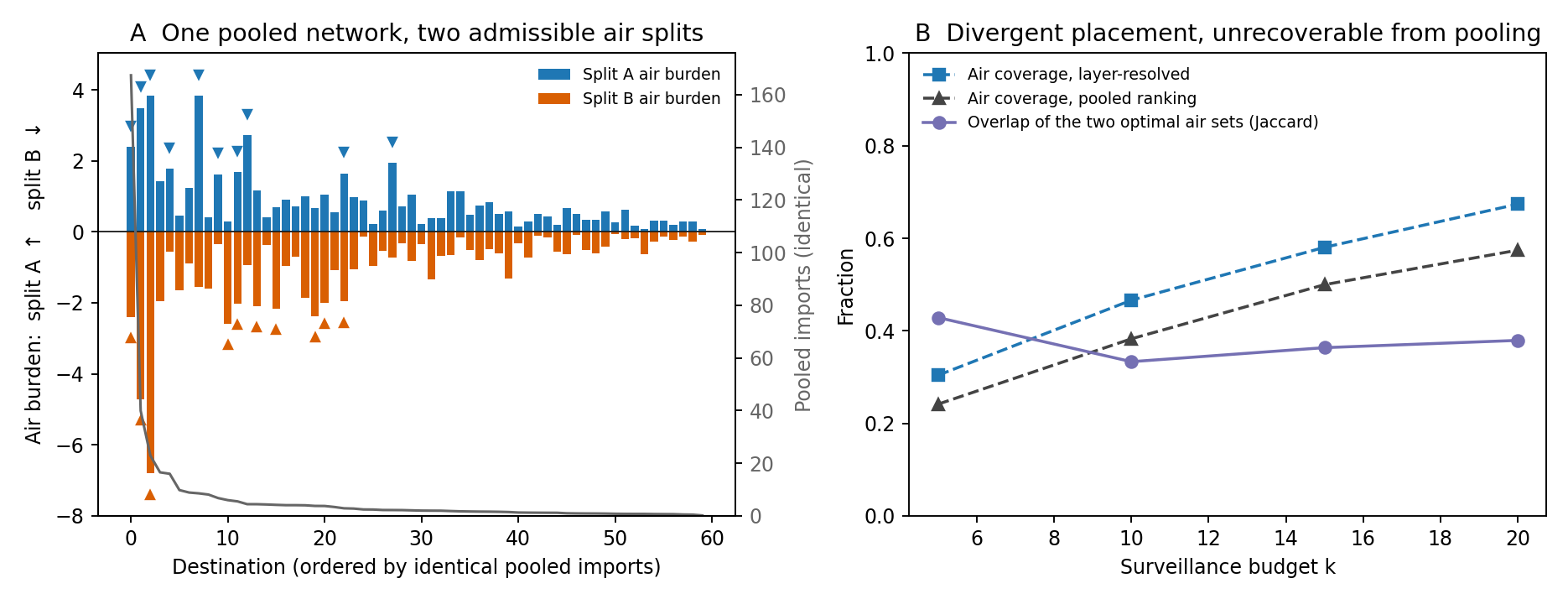}
\caption{Observational-equivalence twin. (A) Two multiplexes share one effective
pooled coupling, so they produce the identical epidemic and identical pooled
per-destination imports (grey line), yet split that coupling into air and
commuting differently (split A upward, split B downward; triangles mark each
split's top-ten air-surveillance targets). (B) Air-import burden of split A
recovered by the pooled and the layer-resolved rankings across budgets $k$, with
the top-$k$ target overlap (Jaccard).}
\label{fig:twin}
\end{figure}

The regret measures what a pooled planner loses against the layer-resolved
truth. When the normalized aggregation map is noninjective, different
route-level compositions are observationally equivalent
(cf. Corollary~\ref{cor:aggregation_spectral}), and we realize this at the epidemic
level with an observational-equivalence twin. Two multiplexes may share the
identical effective pooled coupling, so they generate the identical realized
epidemic and the identical pooled per-destination import counts (maximum
difference zero by construction), yet they split that coupling into air and
commuting differently, the two splits A and B that Figure~\ref{fig:twin}A shows.

Averaged over
200 replicates, the two splits carry a comparable overall air share ($13.1\%$
and $14.9\%$), but they place air on different destinations: their top-ten
air-surveillance targets overlap in only five of ten (Jaccard $0.33$), and the
overlap stays between $0.33$ and $0.43$ across budgets. A planner who sees only
the pooled counts cannot tell the two splits apart, so the pooled ranking
recovers $38.3\%$ of split A's air burden at $k=10$ against $46.6\%$ for the
layer-resolved ranking, a regret consistent with the direct-pooling result.
Figure~\ref{fig:twin}B shows this gap across budgets. The loss is not
a matter of estimator quality. It is the noninjectivity of the aggregation map
made concrete: one pooled observation is consistent with route-level truths that
a mode-specific intervention would treat differently (construction in
the electronic supplementary material, section~S9).

The composite static aggregate gives a deliberately stricter sensitivity
analysis, and it shows the complementary invariant. Its five most-imported
destinations agree with those of the full multiplex, and its top-ten lists
overlap in nine destinations (Jaccard $0.82$). The largest-import destinations
therefore survive pooling, while the route composition collapses to a reweighted
pushforward that pooling cannot recover.

\subsection{Empirical case studies}
\label{sec:real_data}

The synthetic study shows that pooling erases layer attribution whenever the
layers carry distinct route information. Two empirical case studies examine how
the same discrete-time hazard form behaves under different epidemic and data
conditions. In the 2009 A(H1N1) pandemic in the United States, the air layer
records inter-MSA passenger volume during a period of normal domestic service.
In early-2020 Italy, air travel was disrupted and the available air layer is a
2019 service-intensity proxy with incomplete endpoint reporting. The diseases,
observation scales, temporal resolutions, mobility measurements, and sample
sizes differ. These differences prevent a controlled comparison. Each analysis
asks whether a mobility layer adds fitted information beyond epidemic timing and
whether an air-specific ranking retains fitted burden to target once the layers
are pooled. A structural comparison then relates each system's raw multilayer
alignment to resident population, the shared factor of
Theorem~\ref{thm:shared_factor}.

\subsubsection{A shared invasion-hazard model}
\label{sec:hazard_model}

Both reconstructions adapt the ensemble-network perspective of
\citet{zhang2026spatial} to a discrete-time first-invasion hazard. For each
non-seed destination $j$,
\begin{equation}
h_j(t)=b_0+b_tg(t)+c_a\sum_{i\ne j}F^{a}_{ij}\frac{I_i(t)}{N_i}
       +c_c\sum_{i\ne j}F^{c}_{ij}\frac{I_i(t)}{N_i},
\label{eq:invasion_hazard}
\end{equation}
where $F^a$ and $F^c$ are the directed air and commuting flows, $I_i(t)/N_i$ is a
lagged infectious-prevalence proxy at source $i$, and
$g(t)=\sum_i I_i(t)/\sum_i N_i$ is national infectious prevalence, with
$b_0,b_t,c_a,c_c\geq0$. The two background terms carry constant external risk and
the epidemic's national time profile without using network topology. The daily
invasion probability is $1-\exp\{-h_j(t)\}$, and the likelihood is survival
before the observed onset followed by invasion on the onset day. The fitted
coefficients absorb the different units of the two flows, so their magnitudes
have no common physical scale.

We fit five nested specifications. The constant-only baseline retains $b_0$
alone, the timing baseline adds $b_tg(t)$, the air-only and commuting-only models
add a single mobility term, and the full model adds both. Leave-one-region-out
refitting removes the held-out destination from the likelihood and from every
prevalence and background covariate before predicting its onset, which gives a
retrospective target-excluded reconstruction rather than a prospective forecast.
A paired regional bootstrap quantifies the change in absolute onset error. A
source-label permutation null refits after reassigning the origins of the
mobility flows, applying the same permutation to both layers, holding destination
labels fixed, and zeroing self-routes. Its one-sided $p$-value is the fraction of
permuted likelihood gains at least as large as the observed gain. Each
application sets its own onset threshold, prevalence lag, and seeds, specified in
the corresponding application section.

\subsubsection{United States: air-informed reconstruction}
\label{sec:us_app}

We analyze the United States with the metropolitan mobility and incidence data of
\citet{zhang2026spatial}, covering 369 metropolitan statistical areas (MSAs) with
a directed census commuting matrix, a daily inter-MSA air-passenger matrix, and a
weekly per-capita ILI$^{+}$ incidence series over 15 weeks. The air layer records
passenger volume rather than service endpoints. Onset is the first week in which
cumulative incidence reaches a fixed level of the series, with a primary threshold
of $0.002$ and the same conclusions at $0.001$ and $0.005$. The earliest-onset
MSAs are seeds, which leaves 215 non-seed onsets, and the prevalence proxy uses a
one-week lag. The air and commuting flows are almost uncorrelated across region
pairs, with log-flow correlation $0.03$, so each layer can carry independent
introduction information.

At the primary specification, air provides additional fitted information.
Relative to the timing baseline, air alone raises the log-likelihood by $15.8$
units and commuting alone by $19.6$ units, and both
coefficients are positive in the full model. Because the layers are
near-orthogonal, each contributes on top of the other. Air adds $5.2$ units above
the commuting model, and commuting adds $9.0$ above the air model. The
source-label permutation null places the observed air increment at $p=0.017$
against the timing baseline and at $p=0.033$ for the increment above commuting.
Target-excluded onset error falls from $0.70$ weeks under the timing baseline to
$0.61$ weeks under the full model. The route ensemble attributes $11.6\%$ of
post-seed introductions to air (bootstrap $95\%$ interval $8.4$--$15.3\%$),
$15.8\%$ to commuting, and the remainder to external risk (the electronic
supplementary material, section~S6).

The fitted hazard also yields a data-based pooling analysis, the empirical
counterpart of the regret $\mathcal R_k$ in \eqref{eq:regret} with $a_j$ the
fitted air contribution at destination $j$'s observed onset. For each non-seed
MSA we evaluate the air and commuting hazard contributions at its observed onset,
rank MSAs by the air contribution ($\mathcal O_k$) or by the pooled
air-plus-commuting contribution ($\mathcal P_k$), and compare the air burden each
set covers. At $k=10$, $C(\mathcal O_{10})=56.2\%$ and $C(\mathcal P_{10})=45.4\%$,
so $\mathcal R_{10}=10.8$ percentage points, with a $95\%$ interval of
$2.2$--$49.0$ from a region-resampling bootstrap that redraws the 215 onset
targets and refits the hazard on each resample. The two top-ten sets overlap in
six MSAs. The interval is wide and right-skewed because few of the resampled
targets are air-dominated, so the point estimate is directional evidence that the
loss is positive and does not fix its magnitude. The controlled synthetic
analysis of Section~\ref{sec:aggregation_hides} carries the tight quantitative
statement, $\mathcal R_{10}=8.2$ points with interval $7.6$--$8.9$.

\subsubsection{Italy: commuting-informed reconstruction}
\label{sec:italy_app}

We reconstruct early Italian spread from the ISTAT 2011 census commuting matrix
\citep{istat2011commuting} and 2019 OpenSky aviation records
\citep{schaefer2014opensky}. OpenSky reports no passenger counts, so aviation
strength measures service intensity rather than passenger volume, and endpoint
reporting is incomplete. The onset reconstruction uses the directed domestic
matrices aggregated to 21 modern NUTS-2 regions, matched to the EC-JRC incidence
series, while the structural comparison of Section~\ref{sec:structural} uses the
finer 110 NUTS-2010 level-3 units. Record counts, endpoint completeness, and the
NUTS reconciliation are given in the electronic supplementary material, section~S5.

Onset is the first date on which cumulative reported incidence reaches 10 cases
per 100{,}000, between 20 February and 31 May 2020. The two regions whose onsets
fall within three days of the earliest are the initial seeds, which leaves 19
non-seed onsets. Because reported incidence lags infection, the exposure
covariates use a seven-day backward shift, consistent with the roughly five-day
COVID-19 incubation period of \citet{lauer2020incubation}. The prevalence proxy
at source $i$ is the shifted seven-day sum of reported new cases divided by
population.

Route uncertainty is propagated through 2{,}000 realizations. We jitter onset
dates by up to two days, sample coefficient vectors from 1{,}000 region-bootstrap
refits so that boundary coefficients can be positive in some draws, perturb the
lagged prevalence proxy with lognormal reporting noise, and draw a source and
layer for each post-seed invasion in proportion to their hazard contributions.
Route eligibility requires positive source prevalence on that day, matching the
fitted hazard rather than imposing an onset-order restriction. The ensemble is
summarized by route occurrence probabilities, layer shares, and the
destination-specific route entropy $H_j=-\sum_r p_{jr}\log p_{jr}$, taken over the
ensemble occurrence probabilities $p_{jr}$ of the source-layer routes $r$ that
introduce destination $j$. Onset thresholds of 5 and 20 cumulative cases per
100{,}000 provide sensitivity checks.

In the primary Italian specification, commuting contributes fitted information
and the air coefficient is zero. The full directed hazard has negative log
likelihood 51.48, against 56.72 for the timing baseline and 79.32 for the
constant-only baseline. The air-only extension equals the timing baseline, and
the commuting-only extension equals the full model.
Directed commuting therefore adds 5.24 log-likelihood units beyond the national
epidemic curve, and air service adds no maximum-likelihood increment in this
specification. The source-label negative control gives $p=0.0195$ across 2{,}000
permutations for that gain, so the spatial arrangement of commuting sources
carries information beyond epidemic timing (the electronic supplementary material, section~S5).
The zero fitted air increment does not by itself verify that air-mediated spread
was absent. It can reflect disrupted travel, the service-intensity proxy,
incomplete endpoint reporting, the 19 non-seed onsets, or mismatch between
reported onset and infection time.

Target-excluded leave-one-region-out reconstruction gives a mean absolute onset
error of 3.05 days under the full model, 4.74 days under the timing baseline, and
7.53 days under the constant-only baseline. The full-versus-timing reduction is
1.68 days (paired bootstrap $95\%$ interval $0.68$--$2.74$). Across onset
thresholds of 5, 10, and 20 cumulative cases per 100{,}000, the full-model errors
are 2.94, 3.05, and 3.74 days, against 3.83, 4.74, and 5.11 days for the timing
baseline, and lag choices from 0 to 10 days preserve the direction of the
improvement. 

Table~\ref{tab:rolling_validation} reports a rolling-origin validation over 14 expanding folds, beginning from
five observed domestic events and censoring prevalence after each cutoff.
Overall the full-model MAE is 5.14 days, against 5.86 for timing-only, 4.57 for
commuting-only, and 6.43 for air-only. The
full model beats timing-only on 9 of 14 folds with two ties, commuting-only beats
timing-only on 9 of 14 folds, and commuting-only performs best in this sample.
With only 14 folds and a fixed epidemic background after each cutoff, the rolling
comparison is a within-outbreak diagnostic rather than a general forecast
benchmark. A normalized pooled comparator, formed by dividing the air and
commuting exposure surfaces by their training-event means and summing them into
one covariate with a single nonnegative coefficient, removes layer labels but
also constrains the relation between the two surfaces. It has negative log
likelihood 54.50, a 2.22-unit gain over timing-only, with target-excluded MAE
4.11 days and rolling MAE 6.14 days. These values describe that constrained
specification and do not isolate the predictive cost of aggregation.

\begin{table}[h]
\centering\small
\caption{Rolling-origin onset validation. Errors are mean absolute onset errors
across 14 expanding folds. Differences are relative to the timing-only model, and
negative values favor the listed model. The final column gives the number of
folds with an error no greater than seven days.}
\label{tab:rolling_validation}
\begin{tabular}{lrrrr}
\toprule
Model & Folds & MAE (days) & Difference vs timing & Within 7 days\\
\midrule
Timing-only & 14 & 5.86 & 0.00 & 11/14\\
Full air-commuting & 14 & 5.14 & $-0.71$ & 13/14\\
Commuting-only & 14 & 4.57 & $-1.29$ & 14/14\\
Air-only & 14 & 6.43 & $+0.57$ & 10/14\\
Normalized pooled sensitivity & 14 & 6.14 & $+0.29$ & 10/14\\
\bottomrule
\end{tabular}
\end{table}

\begin{figure}[h]
\centering
\includegraphics[width=0.68\textwidth]{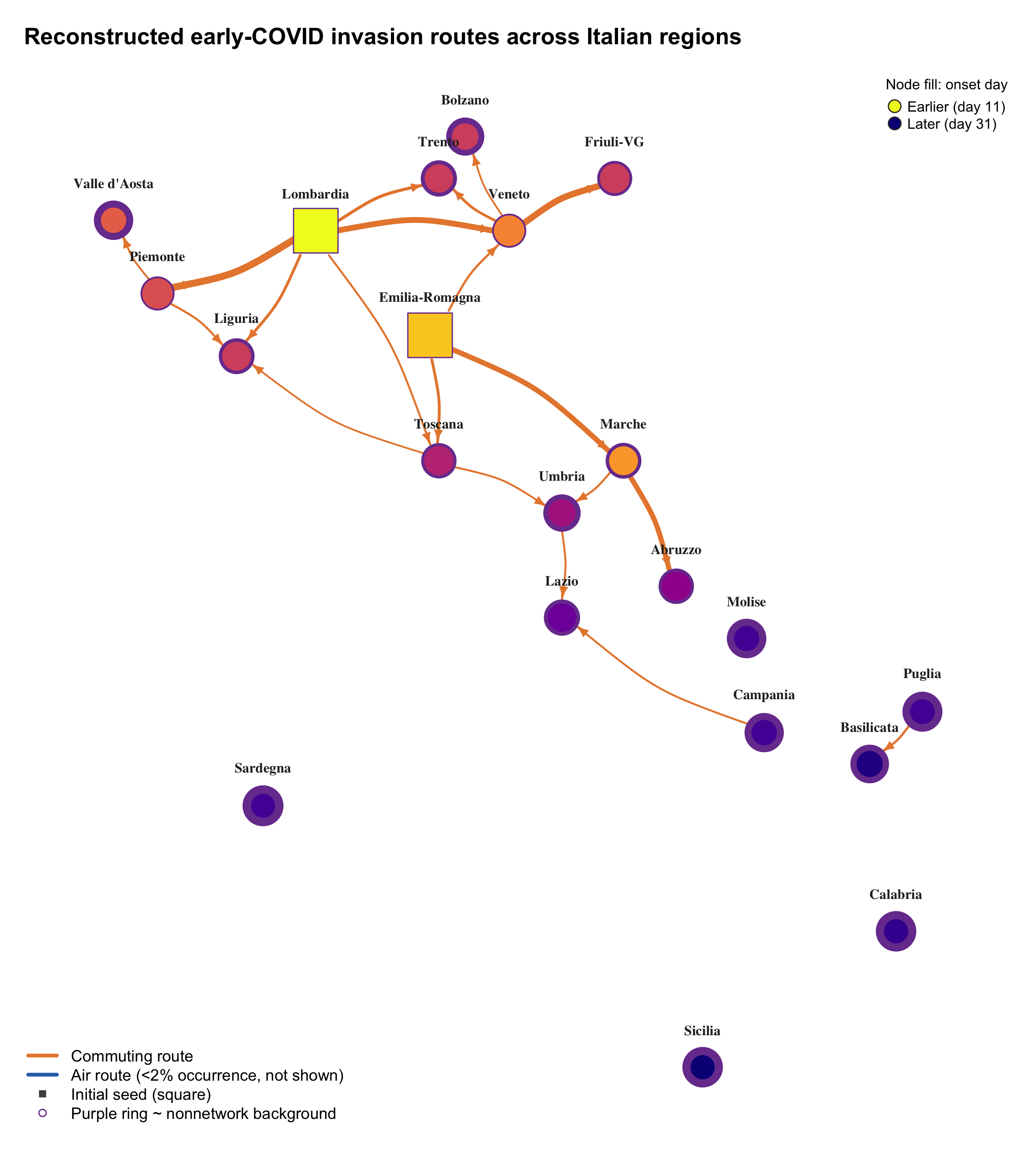}
\caption{Onset-constrained ensemble reconstruction of early Italian invasion,
shown as a directed route network over the 21 NUTS-2 regions. Nodes are at
regional centroids, coloured by onset day, with the two initial seeds (Lombardia
and Emilia-Romagna) drawn as squares. Orange arrows are commuting introductions
with ensemble occurrence at least $0.15$ and width proportional to occurrence
probability; air introductions never exceed $2\%$ occurrence and are omitted.
The purple ring around each region is proportional to its external-introduction
probability.}
\label{fig:route_reconstruction}
\end{figure}

Across the 2{,}000 reconstructed networks, external risk accounts for $52.0\%$ of
post-seed introductions, commuting for $47.5\%$, and air for $0.4\%$ in aggregate.
Figure~\ref{fig:route_reconstruction} makes the spatial pattern explicit. A
commuting backbone links the northern and central regions, and the southern and
island regions are reached mainly by external introductions. At the replicate
level the median shares are $52.6\%$ external, $47.4\%$ commuting, and $0\%$ air,
and the air share has a $95\%$ interval of $0$--$5.3\%$. Layer attribution is
therefore uncertain, and route attribution is more diffuse still. Among 159
domestic routes occurring at least once, $93.1\%$ occur in fewer than $20\%$ of
ensemble members and $0.6\%$ occur in at least $80\%$, and the median
destination-specific route entropy is $0.91$ (the electronic supplementary material, section~S5).

\subsubsection{Comparison of the two cases}
\label{sec:two_regimes}

\begin{figure}[h]
\centering
\includegraphics[width=\textwidth]{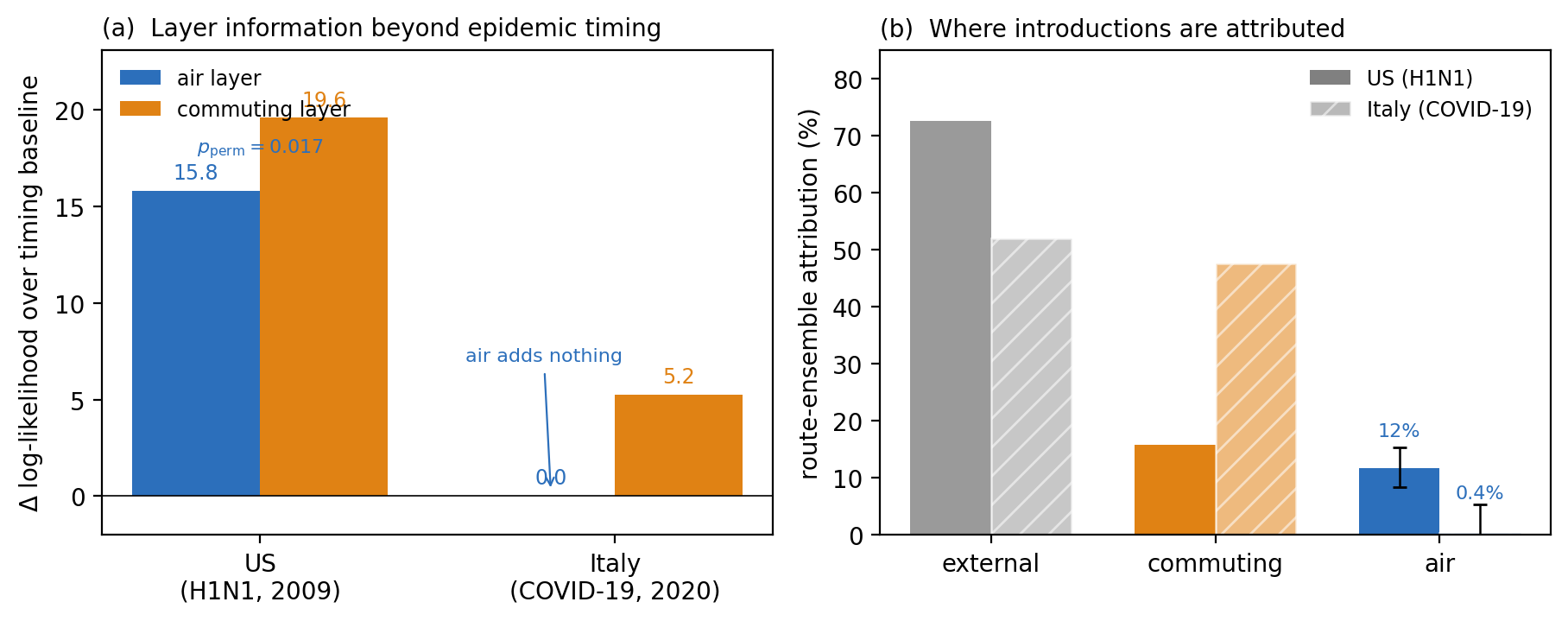}
\caption{Layer attribution in two case studies. (a)~Log-likelihood gain of each
mobility layer over the topology-free epidemic-timing baseline, for the US
A(H1N1) and Italian COVID-19 reconstructions. (b)~Route-ensemble attribution of
post-seed introductions to external risk, commuting, and air, with 95\%
bootstrap intervals on the air share. The US air layer records passenger volume
rather than service endpoints.}
\label{fig:two_regimes}
\end{figure}

Table~\ref{tab:two_regime_audit} and Figure~\ref{fig:two_regimes} place the two
case studies side by side. They use the same hazard form but differ in disease,
measurement, resolution, and sample size. Their comparison illustrates how the
same hazard form behaves across settings; it does not support controlled
cross-setting inference. At the primary US specification, air adds fitted
information beyond commuting, and pooling it into total mobility changes the
destinations selected by the model-based air-surveillance ranking. In Italy the
maximum-likelihood air coefficient is zero, so air has an ensemble introduction
share near zero ($0.4\%$, interval $0$--$5.3\%$) and no comparable fitted
air-targeting diagnostic is defined. These outputs show how layer attribution
depends on the fitted mobility contribution in each case; they do not isolate a
common cross-setting effect of aggregation or establish that air-mediated spread
was absent in Italy.

\begin{table}[h]
\centering\small
\setlength{\tabcolsep}{2pt}
\caption{Fitted-hazard summary of the two case studies, conditional on the
fitted onset models. Air burden covered is the share of fitted air contribution
in the selected destinations, and regret is the layer-resolved minus pooled
coverage $\mathcal R_k$ of \eqref{eq:regret} with a 200-sample region-resampling
bootstrap interval. }
\label{tab:two_regime_audit}
\begin{tabular}{p{0.31\textwidth}p{0.3\textwidth}p{0.3\textwidth}}
\toprule
Quantity & US A(H1N1) & Italy COVID-19\\
\midrule
Non-seed onset targets & 215 & 19\\
Air gain over timing baseline & $+15.8$ & $0.0$\\
Air gain beyond commuting & $+5.2$ & $0.0$\\
Fitted air introduction share & $11.6\%$ [$8.4,15.3$] & $0.4\%$ [$0,5.3$]\\
Air coverage at $k=10$ & $56.2\%$ resolved, $45.4\%$ pooled & No fitted air burden\\
Pooling loss & $10.8$ pp [$2.2,49.0$], $6/10$ overlap & Not applicable, $c_{\rm air}=0$ at MLE\\
\bottomrule
\end{tabular}
\end{table}

\subsubsection{Structural context: population and gateway alignment}
\label{sec:structural}

The reconstructions use the directed flows as inputs. A separate structural
comparison asks where the two layers are jointly prominent in each country,
independent of the epidemic, and connects the raw alignment to the shared-factor
mechanism of Theorem~\ref{thm:shared_factor}. For each layer $\ell$, direction
$A\in\{\mathrm{out},\mathrm{in}\}$, and region $i$, we compute the directed
strength $D_{i\ell}^{A}$ and residualize it on resident population,
\begin{equation}
r_{i\ell}^{A}=\log(1+D_{i\ell}^{A})
-\widehat\alpha_{\ell A}-\widehat\beta_{\ell A}\log N_i,
\label{eq:residual}
\end{equation}
fitting the intercept $\widehat\alpha_{\ell A}$ and slope $\widehat\beta_{\ell A}$
by least squares, with 2009 MSA population for the US and 2019 regional
population for Italy. Raw strengths measure total concentration, and the
residuals isolate functional gateway prominence beyond demographic scale,
proxying $\log G_{i\ell}^{A}$. Inference uses Spearman association, source
permutations that reassign the air residuals while holding commuting fixed,
paired bootstrap intervals, and top-decile overlap against a hypergeometric
independence reference, with Benjamini-Hochberg $q$-values across the six residual
tests in each country. Because zero-service regions force a population-determined
residual, we report an air-active subset in each country, 194 of 369 MSAs and 80
of 110 Italian regions. The Italian commuting and aviation layers are eight years
apart, so the Italian comparison tests persistent regional prominence and makes
no temporal attribution between 2011 commuting and 2019 air traffic.

In the United States, Table~\ref{tab:us_structural} shows a strong raw
association among air-active MSAs, with inbound $\rho=0.563$ and total
$\rho=0.556$. After population adjustment the association reverses sign and
remains significant. The inbound residual correlation is $-0.239$ (permutation
$p=0.0014$, BH $q<0.01$, bootstrap $95\%$ interval $-0.374$ to $-0.094$), the
total residual is $-0.266$, and the outbound residual is $-0.283$. The
full-universe residuals are more strongly negative, from $-0.47$ to $-0.56$, but
that sample mixes served MSAs with the 175 whose zero air strength forces a
population-determined residual, so the air-active estimate is the interpretable
one. A partial-rank correlation given log population reproduces the sign and
magnitude ($-0.20$ air-active, $-0.44$ full universe), and the negative residual
is stable across minimum-passenger thresholds (the electronic supplementary material, section~S6).
Conditional on population, aviation prominence and commuting prominence trade off.
MSAs that carry more air traffic than their size predicts tend to carry less
cross-MSA commuting than their size predicts, consistent with airline-hub metros
and labor-market centers being largely distinct sets of cities. The top-decile
overlap of the two residual rankings is at or below its chance expectation.

\begin{table}[h]
\centering
\small
\caption{Directed air-commuting association across US metropolitan areas. Raw
strengths retain demographic and functional concentration, and
population-residual strengths isolate functional prominence via
\eqref{eq:residual}. The final columns report population-residual top-decile
overlap and its one-sided hypergeometric probability under independent rankings.}
\label{tab:us_structural}
\resizebox{\textwidth}{!}{%
\begin{tabular}{llrrrrrrr}
\toprule
Sample & Direction & $n$ & Raw $\rho$ & Population-residual $\rho$ & Permutation $p$ & BH $q$ & Top-decile overlap & Tail-null $p$\\
\midrule
All MSAs   & Out   & 369 & 0.095 & -0.564 & $<0.001$ & $<0.001$ & 1/37 & 0.984\\
All MSAs   & In    & 369 & 0.201 & -0.466 & $<0.001$ & $<0.001$ & 2/37 & 0.909\\
All MSAs   & Total & 369 & 0.144 & -0.533 & $<0.001$ & $<0.001$ & 1/37 & 0.984\\
Air-active & Out   & 194 & 0.535 & -0.283 & $<0.001$ & $<0.001$ & 2/19 & 0.580\\
Air-active & In    & 194 & 0.563 & -0.239 & 0.0014   & 0.0014   & 3/19 & 0.280\\
Air-active & Total & 194 & 0.556 & -0.266 & $<0.001$ & $<0.001$ & 2/19 & 0.580\\
\bottomrule
\end{tabular}}
\end{table}

Table~\ref{tab:italy} reports the Italian association. The raw air and ground
strengths are positively associated across all 110 regions, especially inbound,
with in-strength correlation $0.466$ and total-strength correlation $0.389$.
After population
adjustment these fall to $0.120$ ($p=0.213$) and $0.064$ ($p=0.501$), and the
outbound residual is $-0.035$ ($p=0.718$), so across the full universe the
evidence for a functional gateway increment is weak. Among the 80 air-active
regions the inbound population-residual association is $0.300$ (permutation
$p=0.0068$, BH $q=0.041$, bootstrap $95\%$ interval $0.065$--$0.499$). Three of
the eight largest inbound commuting residuals are also among the eight largest
inbound air residuals, against $0.8$ expected under independent rankings
(hypergeometric $p=0.030$). Total residual strength is more weakly associated
($\rho=0.222$, $p=0.0438$, BH $q=0.131$), and outbound strength is unrelated
($\rho=0.035$, $p=0.769$). Requiring both flight endpoints to be reported leaves
79 served regions and attenuates the inbound residual to $0.266$ ($p=0.021$,
bootstrap $95\%$ interval $0.039$--$0.464$), with total-strength association
$0.198$ ($p=0.085$). Across thresholds of 10 to 500 annual endpoints the inbound
estimate ranges from $0.176$ to $0.249$. Conditioning on air activity can select
regions through the outcome being studied, so the served-sample result is
consistent with residual gateway alignment among served regions rather than a
demonstration of it. The $n=110$ sample gives no stable estimate of $\lambda_U$
for the full multiplex.

\begin{table}[ht]
\centering
\small
\caption{Directed air-commuting association across Italian NUTS-2010 level-3
regions. Columns are as in Table~\ref{tab:us_structural}: raw strengths retain
demographic and functional concentration, and population-residual strengths
isolate functional prominence via \eqref{eq:residual}.}
\label{tab:italy}
\resizebox{\textwidth}{!}{%
\begin{tabular}{llrrrrrrr}
\toprule
Sample & Direction & $n$ & Raw $\rho$ & Population-residual $\rho$ & Permutation $p$ & BH $q$ & Top-decile overlap & Tail-null $p$\\
\midrule
All regions & Out   & 110 & 0.274 & -0.035 & 0.718 & 0.769 & 1/11 & 0.704\\
All regions & In    & 110 & 0.466 &  0.120 & 0.213 & 0.426 & 3/11 & 0.079\\
All regions & Total & 110 & 0.389 &  0.064 & 0.501 & 0.752 & 2/11 & 0.302\\
Air-active  & Out   &  80 & 0.277 &  0.035 & 0.769 & 0.769 & 2/8  & 0.181\\
Air-active  & In    &  80 & 0.532 &  0.300 & 0.0068& 0.041 & 3/8  & 0.030\\
Air-active  & Total &  80 & 0.443 &  0.222 & 0.0438& 0.131 & 3/8  & 0.030\\
\bottomrule
\end{tabular}}
\end{table}

Population therefore supplies the raw multilayer alignment in both regimes, the
common heavy-tailed factor of Theorem~\ref{thm:shared_factor}. The functional
gateway increment that remains after adjustment differs in sign, negative in the
United States and weakly positive among air-served Italian regions, yet the
message for aggregation is the same. The two layers encode different, and in the
US even opposed, node-prominence information, so pooling them discards that
distinction.

\section{Discussion}

The main result is a coarse-graining principle for heterogeneous multilayer
networks. A shared heavy-tailed node factor generates an explicit interlayer tail
coefficient, and a random diagonal map with Poissonization preserves the radial
tail index while transforming the angular measure. Nonnegative route aggregation
preserves that index but retains only the norm-weighted pushforward in
\eqref{eq:aggregate_spectral}. Invariance of extreme scale and identifiability of
extreme composition are therefore separate properties.

When the normalized aggregation map $T_S$ is noninjective, distinct route-level
spectral measures produce the same aggregate, and no estimator on the aggregate
alone can separate them. Full route-layer recovery is possible exactly when the
retained route signatures are affinely independent
(Proposition~\ref{prop:identifiability}); the route-level measure is then
recoverable in closed form. Aggregation still identifies any target represented
by an aggregate coordinate, such as total importation into a destination.
Equation~\eqref{eq:od_pushforward} localizes the loss under layer pooling: the
source-destination spectral measure remains identified, while the mode split
within each pair does not. Collapsing the multiplex to one coordinate destroys
all route composition. The
first-invasion theorem adds a dynamical reading, in which route shares are
reweighted by establishment probabilities and the radial extreme sets the
accelerated clock.

The experiments realize these distinctions in finite networks. Mixed-Poisson
sampling and nonnegative aggregation leave radial-index estimates near their
generating values, while two allocation mechanisms with identical source
aggregates give markedly different dominant-route distributions. In the
metapopulation model, pooling preserves the regional import-risk ranking but
removes the dominant transport layer for $23\%$ of destinations. The
invasion-acceleration effect is not universal, being strong for hub seeds and
weak or reversed for median and random seeds, so the reliable consequence of
multilayer structure is attribution rather than automatic acceleration.

The controlled experiment sets the within-model size of the pooling cost, and
the two applications show where that cost appears in observational data. At the
primary US specification, air adds information beyond commuting, and the
layer-resolved air ranking covers more fitted air burden than the pooled one, so
pooling misdirects a mode-specific air-surveillance ranking. In the early Italian
wave the fitted air coefficient is zero while lagged commuting improves
target-excluded onset reconstruction. Under disrupted travel, a service-intensity
air proxy, and incomplete endpoints, air attribution is not recoverable there, an
empirical instance of the nonidentifiability the theory describes rather than a
demonstration that air played no role. The structural comparisons add node-level
context. A common demographic factor produces the raw air-commuting alignment in
both countries, matching the shared-factor mechanism, while the functional
residual is weakly positive among air-served Italian regions and negative across
US metropolitan areas. Either way the raw alignment is largely a population
artifact, so pooling discards attribution rather than redundant structure.

For surveillance design the distinction is operational. When the planning target
is a total, the burden imported into a region, its hotspot ranking, or its
arrival time, that quantity is an aggregate coordinate and pooled mobility
suffices. When the target is mode-specific, where to place air-focused screening,
port-of-entry checks, or a layer-targeted control, the relevant quantity is the
angular composition that pooling destroys, and layer-resolved flows are needed.
The regret $\mathcal R_k$ measures the price of ignoring this distinction: a
mode-specific air resource sited from pooled data forfeits about a fifth of the
attainable air coverage in the synthetic system, and it misses air-dominated
destinations because the largest total-import sites are commuting hubs. These are
planning quantities under a stated model, not causal intervention effects.

The theory is finite-dimensional and asymptotic, and its transfer results assume
independence and moment conditions that warrant checking in a given system. The
route-allocation and aggregation matrices are fixed, and letting them coevolve
with the heavy-tailed state is a natural extension. The metapopulation ranking
result is conditional on the simulated data-generating process, and the
aggregation theorem concerns the tail index rather than a finite-sample rank. The
two applications differ in disease, observation scale, temporal resolution,
mobility measurement, and sample size, so they illustrate the same hazard form in
contrasting regimes rather than forming a controlled comparison or estimating a
common cross-setting effect. Their two mobility layers differ in date and units,
air endpoints proxy passenger volume imperfectly, and onset data do not directly
observe infectious movement. The fitted-hazard analysis ranks observed-onset
destinations under one model, so its empirical route probabilities are directional
model outputs rather than precise or causal effect estimates, and they neither
validate an intervention nor identify unobserved transmission links. These limits
affect calibration and not the aggregation result. Coarse-graining can conserve
the scale of an extreme event while destroying its direction, and whenever
mechanism, channel, or path matters, layer identity is part of the state rather
than dispensable metadata.

\appendix
\section*{Ethics}
This work used only publicly available, aggregated, and de-identified mobility
and epidemiological data together with numerical simulations. It did not involve
human participants, identifiable personal data, or animal experimentation, so no
ethical approval was required.

\section*{Data accessibility}
Italian ISTAT 2011 commuting matrices are publicly available from the
\href{https://www.istat.it/non-categorizzato/matrici-del-pendolarismo/}{ISTAT
commuting archive}. OpenSky flight records for 2019 are
available from the \href{https://doi.org/10.5281/zenodo.5815448}{Zenodo
OpenSky archive} and related records. Italian COVID-19 NUTS case series are
available from the \href{https://github.com/ec-jrc/COVID-19}{EC-JRC COVID-19
repository}. Eurostat NUTS-3 population data are available from
the \href{https://ec.europa.eu/eurostat/data/database}{Eurostat database}. The
US metropolitan commuting, air-passenger, and ILI$^{+}$ incidence data used in
the US application are the public release of \citet{zhang2026spatial}. 

Relevant codes and numerical outputs supporting the figures
and reported results are supplied in the
GitLab repository, 
\href{https://gitlab.com/tw3981/multilayer-epidemics}{\texttt{https://gitlab.com/tw3981/multilayer-epidemics}}. 
Third-party source data are not redistributed in the
supplementary package or project repository and must be obtained from the
public providers cited above under their respective terms.

\section*{Declaration of AI use}
The authors used an AI-assisted tool (Claude Opus 4.8, Anthropic) for language editing
and for help drafting and refining prose, LaTeX, and analysis code. The authors
designed the study, developed the theory and proofs, specified and ran the
simulations and empirical analyses, and verified every result. No AI system was
used to generate scientific insights, analyze data, or draw conclusions, and no
AI system is listed as an author. The authors take full responsibility for the
content and correctness of the manuscript.

\section*{Conflict of interest}
The authors declare no competing interests.

\section*{Authors' contributions}
\textbf{Tiandong Wang:} Conceptualization, Methodology, Software, Formal
analysis, Data curation, Visualization, Writing--original draft.
\textbf{Wei Yang:} Conceptualization, Methodology, Supervision,
Writing--review \& editing.

\section*{Funding}
This work was supported by National Key R\&D Program of China under Grant No.
2025YFA1016503. T. Wang was also funded by the National Natural Science Foundation
of China under Grant 12301660 and the Science and Technology Commission of Shanghai
Municipality under Grant 23JC1400700.

\bibliographystyle{plainnat}
\bibliography{references}

\end{document}


\maketitle

This document gives the finite-sample designs, seed-stratified epidemic
results, synthetic temporal and population-tail sensitivity analyses, and the
Italian and US application details summarized in the manuscript. The required
contents and verification steps for the public reproducibility release are
listed in \texttt{REPRODUCIBILITY\_CHECKLIST.md}.

\section*{S1. Finite-sample transfer experiments}

The mixed-Poisson experiment generated 600 networks of 110 nodes. A Pareto
source factor with tail index $\tau=1.6$ was multiplied by lognormal gateway
factors and sampled with independent Poisson noise. Hill estimates were
computed at thresholds from 2\% to 20\%; the median estimate at the prespecified
10\% threshold was 1.61 (simulation interval 1.50--1.75). Layer-specific
Gamma multipliers reduced cross-layer top-decile overlap while preserving the
radial index. Population residualization used a prespecified gateway increment;
the null rejection rate was 0.050 and power was 0.946 in the simulated panels.
Short log-Gaussian panels underestimated the residual tail coefficient when
the true dependence was strong, motivating the common-mode and block-bootstrap
checks in the main analysis.

\section*{S2. Route allocation and aggregation}

For each of 800 heavy-tailed source--layer totals, we allocated mass over
directed routes using either a concentrated gateway allocation or a diffuse
allocation with the same source aggregates. Mixed-Poisson route counts were
then aggregated by source and layer. Across thresholds, the recovered radial
index remained near 1.6. The dominant-route distributions differed despite
identical source totals (total-variation distance 0.566), demonstrating the
nonidentifiability of route composition after aggregation.

\begin{figure}[ht]
\centering
\includegraphics[width=0.9\textwidth]{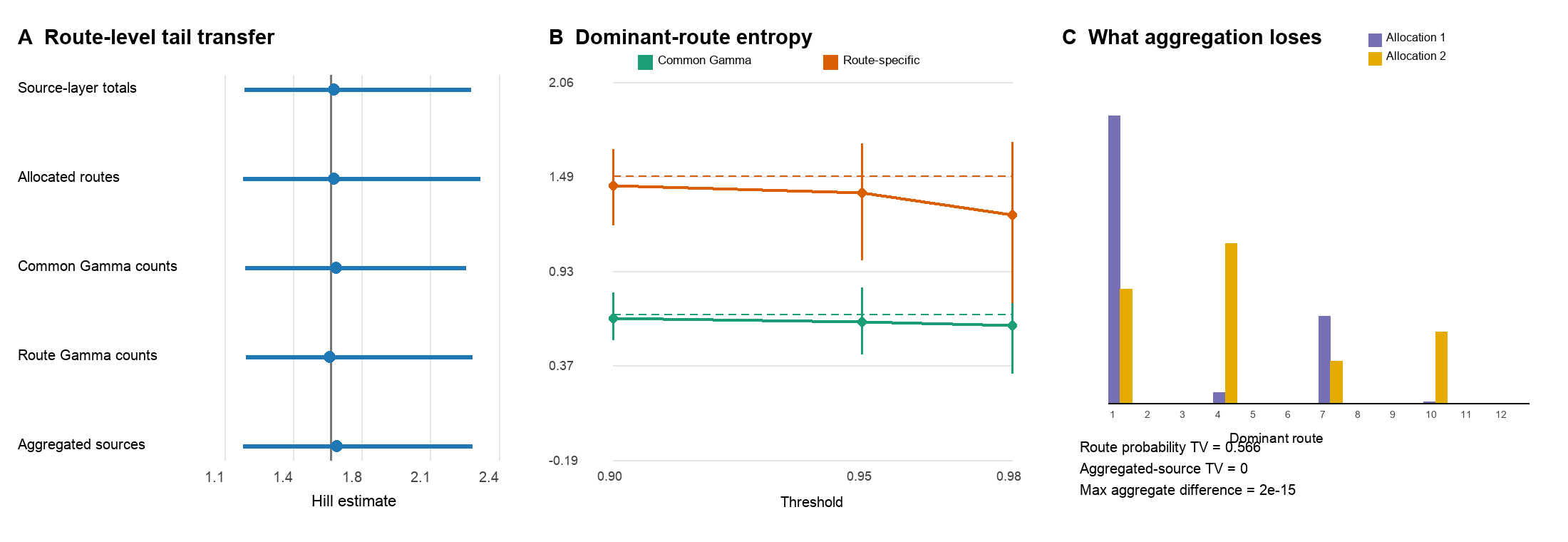}
\caption{Finite-sample route allocation and aggregation checks. The two
allocation mechanisms have identical source aggregates but different
dominant-route distributions.}
\end{figure}

\section*{S3. Seed-placement sensitivity}

The epidemic experiment used 12 independently generated 60-region networks
and 20 paired replicates per network. Each network was simulated with five
seed strategies: aviation hub, commuting hub, population hub, median-population
region, and random region. Table~\ref{tab:seed} reports the median and
interquartile range (IQR) across network means for the full-minus-modified
difference. Negative timing values indicate earlier secondary invasion in the
full multiplex; positive final-size values indicate more invaded regions.

\begin{table}[ht]
\centering\small
\caption{Seed-stratified epidemic differences.}
\label{tab:seed}
\resizebox{\textwidth}{!}{\begin{tabular}{llrrr}
\toprule
Seed & Removed alignment & Timing median [IQR] & Final median [IQR] & $\Pr(\text{timing}<0)$\\
\midrule
Aviation hub & Population & $-17.7$ [$-21.5,-8.5$] & $8.4$ [$4.8,12.1$] & 1.00\\
Aviation hub & Gateway & $0.9$ [$-7.0,3.6$] & $-0.1$ [$-1.8,2.4$] & 0.42\\
Commuting hub & Population & $-9.9$ [$-19.0,-9.3$] & $6.2$ [$5.2,7.0$] & 1.00\\
Commuting hub & Gateway & $-11.4$ [$-19.0,-4.7$] & $7.3$ [$3.6,9.1$] & 0.92\\
Population hub & Population & $-20.1$ [$-22.3,-18.5$] & $8.1$ [$5.6,10.5$] & 1.00\\
Population hub & Gateway & $3.3$ [$-4.3,4.2$] & $-3.0$ [$-5.3,-1.3$] & 0.33\\
Median population & Population & $-0.1$ [$-3.0,1.9$] & $-1.9$ [$-2.8,-1.5$] & 0.50\\
Median population & Gateway & $3.4$ [$-6.0,5.3$] & $-2.6$ [$-4.6,-0.4$] & 0.42\\
Random region & Population & $-2.9$ [$-6.3,2.7$] & $-3.4$ [$-5.0,-0.6$] & 0.58\\
Random region & Gateway & $0.1$ [$-6.7,1.1$] & $-1.4$ [$-3.5,1.0$] & 0.50\\
\bottomrule
\end{tabular}}
\end{table}

\section*{S4. Synthetic temporal sensitivity}

This section reports the sensitivity of the synthetic epidemic to the temporal
spatial fields. Each layer's daily multiplicative field decomposes into a
spatially constant common mode and an independent spatial residual; the
common-mode fraction $f_\ell$ is the share of log-field variance carried by the
constant term, so $(1-f_\ell)/2$ measures the residual spatial contrast.

The two temporal spatial comparisons in the synthetic epidemic have small,
short-lived effects. They advance first secondary invasion by about $1.7$ days
and produce $0.36$--$0.43$ more regions by day 60. Day-110 differences stay
within $0.31$ region, with intervals spanning zero. In this calibration,
observed spatial placement slows early spread slightly. The sign of the
epidemic effect depends on local growth and persistent structural heterogeneity.
The weak effect occurs under
$(f_{\mathrm{air}}=0.40,f_{\mathrm{commute}}=0.25)$, which leaves spatial
contrast $(1-f)/2$ of order $0.30$--$0.38$. Demographic heterogeneity (Pareto
index $1.8$) generates large persistent hub differentials. The temporal spatial
channel would become more consequential when demographic heterogeneity is
compressed, the epidemic seed lies in a mid-sized hub, and the spatial range
approaches the characteristic layer reach.

In a high-common-mode scenario ($f\to0.9$), the
spatial contrast $(1-f)/2$ falls to $0.05$. Removing temporal variation alone
advances secondary invasion by $1.25$ days, but changes final invasion by only
$0.11$ region (95\% interval $-0.48$ to $0.70$). Setting the layer-specific
transmission efficiencies to a common value accounts for a $3.56$-region final difference
($3.00$--$4.13$), close to the composite static aggregate's $3.65$-region
difference ($3.05$--$4.24$). The composite also advances first-secondary
timing by $2.19$ days. The long-run aggregate difference is therefore driven
mainly by using a common transmission efficiency across layers, not by temporal
averaging or loss of route labels alone. These effect sizes describe the
specified synthetic mechanism experiment and are not calibrated estimates for
Italy or COVID-19.

\section*{S5. Italian lag, target exclusion, and coefficient uncertainty}

\paragraph{Data provenance.} The commuting layer comes from the ISTAT 2011 census
origin-destination matrix of work journeys \citep{istat2011commuting}, with
residence and usual workplace defining the directed origin and destination.
Municipality flows are aggregated to the 110 Italian NUTS-2010 level-3 units,
which match the 2011 province system \citep{eurostatnuts2010} and are finer than
the 21 NUTS-2 regions used for the onset reconstruction. The layer contains
19{,}172{,}014 workers on 3{,}988 directed regional pairs including self-flows,
and the structural association uses the 1{,}744{,}083 workers who cross a NUTS-3
boundary, excluding within-region commuting. The aviation layer comes from 2019
OpenSky records \citep{schaefer2014opensky}. The scan of 31{,}218{,}164 global
records identifies 121 active Italian airports and, after mapping endpoints to
NUTS-3 regions, 540{,}548 departures and 594{,}552 arrivals associated with
Italy, international counterparts included. The domestic matrix used in the
reconstruction contains 138{,}314 flights on 1{,}684 directed regional pairs.
Endpoint reporting is incomplete. In total 115{,}863 Italian departures lack a
reported destination and 154{,}397 arrivals lack a reported origin, and a
complete-endpoint sensitivity analysis retains the 424{,}685 departures and
440{,}155 arrivals for which both endpoints are reported. A region is air-active
if it has at least one observed 2019 endpoint, which gives 80 of the 110 regions.

Reported incidence was shifted backward by the specified lag before forming
source prevalence. The primary seven-day lag was motivated by the approximately
five-day incubation estimate of \citet{lauer2020incubation} and is not treated
as a fitted parameter. Table~\ref{tab:lag} gives target-excluded leave-one-region-out
errors across lag choices. Each held-out destination was removed from the
background denominator and all exposure covariates before prediction.

\begin{table}[ht]
\centering\small
\caption{Lag sensitivity of target-excluded onset reconstruction.}
\label{tab:lag}
\begin{tabular}{rrrr}
\toprule
Lag (days) & Full MAE & Timing MAE & Reduction\\
\midrule
0 & 3.37 & 4.68 & 1.32\\
3 & 3.37 & 4.32 & 0.95\\
5 & 3.16 & 4.47 & 1.32\\
7 & 3.05 & 4.74 & 1.68\\
10 & 2.74 & 4.47 & 1.74\\
\bottomrule
\end{tabular}
\end{table}

The seven-day primary fit has negative log likelihood 51.48 (timing baseline
56.72), and the target-excluded MAE reduction is 1.68 days (bootstrap 95\%
interval 0.68--2.74). The source-label negative control uses 2,000 paired
permutations ($p=0.0195$). The constrained maximum-likelihood air coefficient
is zero, but region-bootstrap refits place it above zero in 10.6\% of draws
(median 0; 95\% interval 0--0.0131). This active-set uncertainty is propagated
through 2,000 route-ensemble realizations. Aggregate ensemble attribution is
52.0\% external, 47.5\% commuting, and 0.4\% air; the replicate-level air
share has a 0--5.3\% interval.

\begin{figure}[ht]
\centering
\includegraphics[width=\textwidth]{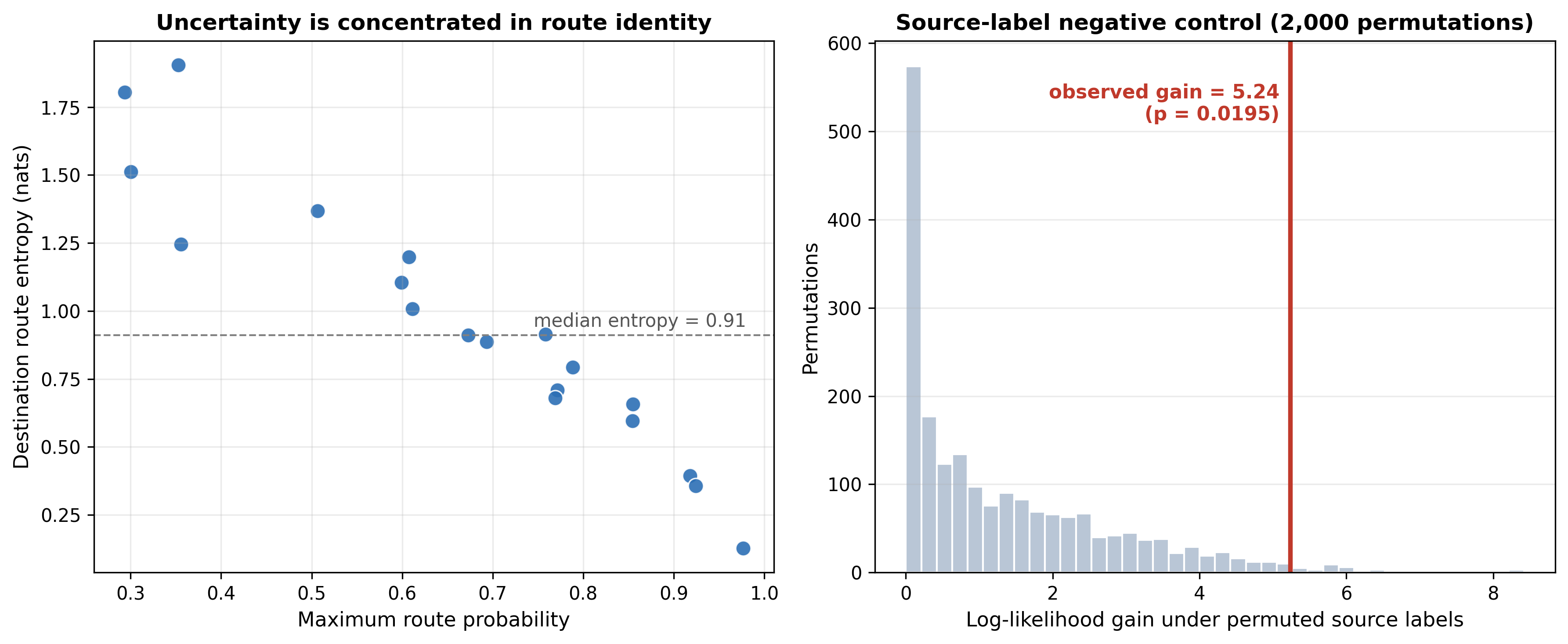}
\caption{Route uncertainty in the Italian onset reconstruction. Left:
destination-specific route entropy versus the maximum route probability across
the 2,000-member ensemble; the dashed line marks the median entropy. Right:
source-label permutation null of the log-likelihood gain of the full mobility
model over the topology-free epidemic-timing baseline, with the observed gain
(red) at $p=0.0195$.}
\label{fig:si_route_uncertainty}
\end{figure}

Figure~\ref{fig:si_route_uncertainty} summarizes the resulting route
uncertainty. Destination-specific route entropy decreases with the maximum
route probability, so the uncertainty is concentrated in the identity of the
introducing route rather than in whether an introduction occurs; the median
destination entropy is $0.91$. The source-label negative control shows that the
observed $5.24$-unit log-likelihood gain of the full mobility model over the
topology-free epidemic-timing baseline is exceeded by only $38$ of $2{,}000$
permuted reassignments ($p=0.0195$).

\section*{S6. US pandemic-influenza application: air-informed onset reconstruction}

This section details the US application summarized in the main text. The
onset-reconstruction estimators assign a positive fitted contribution to air at
the primary US specification. Italy provides a contrasting case with a zero
fitted air coefficient. Differences in disease, measurement, resolution, and
sample size prevent a controlled comparison.

\paragraph{Data.} We use the public US metropolitan mobility and incidence
release of \citet{zhang2026spatial}: a directed $369\times369$ census commuting
matrix (commuters resident in MSA $i$ working in MSA $j$), a daily inter-MSA
air-\emph{passenger} matrix (aggregated to a mean daily flow), a weekly
per-capita ILI$^{+}$ incidence series over 15 weeks during the 2009 A(H1N1)
pandemic, and MSA resident populations. The air layer is passenger volume, not
service intensity, so this control also addresses the OpenSky proxy limitation
of the European data.

\paragraph{Model and estimators.} We fit the same discrete-time invasion hazard
$h_j(t)=b_0+b_tg(t)+c_a\sum_{i\ne j}F^{a}_{ij}I_i(t)/N_i
+c_c\sum_{i\ne j}F^{c}_{ij}I_i(t)/N_i$
as for Italy in the main text, with prevalence proxy $I_i(t)/N_i$ given by
the incidence series at a one-week backward lag and national timing covariate
$g(t)$ equal to the population-weighted mean incidence. Flows enter as
$\sum_{i\ne j}F^{a}_{ij}I_i(t)/N_i$ and $\sum_{i\ne j}F^{c}_{ij}I_i(t)/N_i$,
with each flow matrix scaled by its nonzero median so coefficients are of order
one. Onset is the first week in which cumulative incidence reaches a fixed
level; the earliest-onset MSAs are seeds. At the primary threshold $0.002$ there
are 5 seeds and 215 non-seed onset targets. We report nested log-likelihood
comparisons (constant, timing, air, commuting, full), target-excluded
leave-one-region-out onset error, a source-label permutation control (1{,}000
reassignments of flow origins), and a route ensemble (2{,}000 realizations
attributing each post-seed onset to external, air, or commuting risk in
proportion to the fitted hazard contributions).

\paragraph{Results.} The air and commuting layers are near-orthogonal across
region pairs (log-flow correlation $0.03$), so each can carry independent
introduction signal. At the primary specification, both improve fit.
Table~\ref{tab:us_control} reports the nested
fit across onset thresholds. At the primary threshold, air adds $15.8$
log-likelihood units over the timing baseline and commuting $19.6$; in the full
model both coefficients are positive, air adds $5.2$ units beyond commuting, and
commuting $9.0$ units beyond air. The source-label permutation null is exceeded
at $p=0.017$ for the air-versus-timing increment and $p=0.033$ for the
air-beyond-commuting increment. Target-excluded onset error falls from $0.70$ to
$0.61$ weeks. The route ensemble attributes $11.6\%$ of introductions to air
(interval $8.4$--$15.3\%$), $15.8\%$ to commuting, and $72.6\%$ to external
risk. The air log-likelihood increment over timing is positive at every
threshold ($+7.7$, $+15.8$, $+10.8$ at $0.001$, $0.002$, $0.005$); the
source-label significance of the air-only increment is strongest at the primary
threshold, where the onset times are best spread across the 15-week window.
Because incidence is weekly, onset resolution is coarse and the absolute onset
errors are correspondingly small; the log-likelihood increments and the
attribution shares are the more sensitive summaries.

\paragraph{Pooling analysis.} We translate the fitted layer contributions into a
surveillance-targeting diagnostic. For each non-seed MSA at its observed onset,
we calculate the fitted air contribution and the fitted air-plus-commuting
contribution. A layer-resolved air-surveillance ranking orders MSAs by the first
quantity, while a pooled ranking orders them by the second. At a budget of ten
MSAs, the layer-resolved ranking captures $56.2\%$ of fitted air burden and the
pooled ranking captures $45.4\%$, for a $10.8$ percentage-point loss. The sets
share six MSAs. We refit the full hazard in 200 regional bootstrap samples and
apply each coefficient draw to the original targets; the corresponding 95\%
interval for the loss is $2.2$--$49.0$ percentage points. Table~\ref{tab:us_pool}
reports budgets of 5, 10, and 20 MSAs. This analysis uses modelled hazard
contributions and does not estimate an intervention effect.

\begin{table}[ht]
\centering\small
\caption{US A(H1N1) onset-model fit across thresholds. Log-likelihood
gains are relative to the topology-free timing baseline; positive values favor
the listed model. ``Air $\mid$ commuting'' is the additional gain of air on top
of the commuting-only model. MAE is target-excluded leave-one-region-out onset
error in weeks.}
\label{tab:us_control}
\begin{tabular}{lrrr}
\toprule
Quantity & Thr.\ $0.001$ & Thr.\ $0.002$ & Thr.\ $0.005$\\
\midrule
Non-seed onset targets & 206 & 215 & 219\\
Air gain over timing & $+7.7$ & $+15.8$ & $+10.8$\\
Commuting gain over timing & $+8.6$ & $+19.6$ & $+25.5$\\
Full gain over timing & $+11.4$ & $+24.8$ & $+26.3$\\
Air $\mid$ commuting gain & n/a & $+5.2$ & n/a\\
Air source-label $p$ (vs timing) & $0.37$ & $0.017$ & $0.34$\\
MAE, timing / full (weeks) & $0.73/0.67$ & $0.70/0.61$ & $0.81/0.64$\\
\bottomrule
\end{tabular}
\end{table}

\begin{table}[ht]
\centering\small
\caption{US H1N1 fitted-hazard pooling analysis. Air burden covered is the share
of fitted air contribution in the selected destinations. Regret is the
layer-resolved minus pooled coverage; intervals come from 200 region-resampling
bootstrap refits.}
\label{tab:us_pool}
\begin{tabular}{rrrrr}
\toprule
$k$ & Layer-resolved (\%) & Pooled (\%) & Regret (pp) & Top-$k$ overlap\\
\midrule
5  & 40.7 & 30.6 & 10.1 [0.0, 35.9] & 2/5\\
10 & 56.2 & 45.4 & 10.8 [2.2, 49.0] & 6/10\\
20 & 70.4 & 62.5 & 8.0 [1.6, 34.8] & 12/20\\
\bottomrule
\end{tabular}
\end{table}

The Italian counterparts of these quantities are an air gain over timing of
zero, an air source-label increment that cannot be estimated because the
maximum-likelihood air coefficient is zero, and an ensemble air share of $0.4\%$
(interval $0$--$5.3\%$). These heterogeneous applications illustrate fitted
layer attribution and do not establish a common cross-setting effect. The
Italian zero estimate can also reflect disrupted travel, the service-intensity
proxy, incomplete endpoints, the small sample, or onset-time measurement.

\paragraph{Structural gateway association.} The main text reports a directed
air--commuting structural comparison for the US metropolitan areas, the
counterpart to the Italian association analysis, using the same population
residual $r_{i\ell}^{A}=\log(1+D_{i\ell}^{A})-\widehat\alpha_{\ell A}
-\widehat\beta_{\ell A}\log N_i$. This note records its robustness. Raw air and
commuting strengths are strongly aligned among the 194 air-active MSAs (total
Spearman $0.556$), but the alignment is demographic: after residualizing on log
population the association turns negative ($-0.266$ total, $-0.239$ inbound;
permutation $p<0.01$). An independent partial-rank correlation given log
population reproduces the sign and magnitude ($-0.20$ among air-active MSAs,
$-0.44$ across all 369). The negative residual is stable across minimum-service
thresholds: restricting to MSAs above $100$ to $2{,}000$ mean daily air
passengers leaves the inbound residual between $-0.22$ and $-0.30$ and the total
residual between $-0.23$ and $-0.31$, all with permutation $p<0.05$ ($n$ from
$186$ to $92$). The full-universe residuals are more strongly negative because
the 175 MSAs without scheduled service carry a population-determined air
residual; the air-active estimate is therefore the interpretable one. Because the
US air layer is passenger volume, this measures the functional gateway increment
more directly than the OpenSky service proxy allows for Italy, and its sign shows
that aviation prominence and commuting prominence occupy largely distinct
metropolitan areas once population is removed.

\section*{S7. Decision cost of pooling in the synthetic epidemic}

This section details the surveillance-regret calculation summarized in the main
text. It quantifies, within the Full-multiplex simulation, how much of the
air-import burden a mode-specific surveillance allocation forfeits when it must
be targeted from pooled rather than layer-resolved mobility. It is a
within-model classification exercise and makes no causal claim.

\paragraph{Setup.} For each destination we aggregate infectious imports by layer
over the 300 Full-multiplex replicates, giving expected air and commuting import
counts. A planner allocates a mode-specific air resource to $k$ destinations. The
\emph{layer-resolved} (oracle) planner ranks destinations by expected air-import
burden; the \emph{pooled} planner cannot observe the air layer and ranks by
expected total imports. Writing $a_j$ for destination $j$'s expected air-import
burden and $t_j=a_j+c_j$ for its expected total imports, let $\mathcal O_k$ and
$\mathcal P_k$ be the top-$k$ destinations by $a_j$ and by $t_j$, and let
$C(\mathcal S)=\sum_{j\in\mathcal S}a_j\big/\sum_j a_j$ be the fraction of network
air-import burden a set $\mathcal S$ covers. The regret is the air coverage the
pooled ranking forgoes,
\begin{equation}
\mathcal R_k=C(\mathcal O_k)-C(\mathcal P_k).
\label{eq:si_regret}
\end{equation}
We report $C(\mathcal O_k)$, $C(\mathcal P_k)$, and $\mathcal R_k$, the overlap of
the two top-$k$ sets, and how many of the pooled planner's picks are
air-dominated ($\mathbb E[\text{air}]>\mathbb E[\text{commuting}]$). Regret
intervals are paired bootstrap over replicates ($2{,}000$ resamples).

\paragraph{Result.} The simulation has 60 destinations, of which 14 ($23\%$) are
air-dominated, with a $14\%$ network-wide air share that matches the aggregation
analysis in the main text. Table~\ref{tab:regret} reports the sweep. Pooling
forfeits $7$--$9$ percentage points of attainable air coverage across budgets,
and the pooled top-$k$ contains essentially none of the air-dominated
destinations, because the largest total-import destinations are commuting hubs.
The loss is therefore structural, not a matter of a few swapped ranks: the
pooled resource is placed where air is not the operative pathway.

\begin{table}[ht]
\centering\small
\caption{Surveillance regret of pooled versus layer-resolved targeting in the
synthetic epidemic. ``Air covered'' is the fraction of network air-import burden
in the top-$k$ destinations under each ranking; regret is their difference in
percentage points with a paired bootstrap 95\% interval. ``Overlap'' is the size
of the intersection of the two top-$k$ sets. ``Air-dominated in pooled top-$k$''
counts how many of the pooled picks have air as their dominant import layer, out
of 14 air-dominated destinations network-wide.}
\label{tab:regret}
\begin{tabular}{rrrrrr}
\toprule
$k$ & Oracle air covered & Pooled air covered & Regret (pp) & Top-$k$ overlap & Air-dom.\ in pooled top-$k$\\
\midrule
5  & $30.8\%$ & $24.0\%$ & $6.8$ $[6.2,7.4]$ & $3/5$   & $0$ of $14$\\
10 & $46.5\%$ & $38.3\%$ & $8.2$ $[7.6,8.9]$ & $6/10$  & $0$ of $14$\\
15 & $57.8\%$ & $50.7\%$ & $7.0$ $[6.5,8.5]$ & $11/15$ & $1$ of $14$\\
20 & $67.1\%$ & $57.9\%$ & $9.2$ $[8.4,10.0]$& $13/20$ & $1$ of $14$\\
\bottomrule
\end{tabular}
\end{table}

\section*{S8. Numerical validation of the first successful-event theorem}

This section validates the first successful-event route and time limits of the
main text. The theorem predicts that, conditional on the total route intensity
being extreme, the first-success route distribution converges to
$\pi^*_e=\int_{\mathbb S_+}(s_ew_e/\bm s^\top\bm w)\,H_R(d\bm w)$, the rescaled
first-success time survival to
$S^*(x)=\int_{\mathbb S_+}\int_1^\infty e^{-xr\,\bm s^\top\bm w}\tau
r^{-\tau-1}\,dr\,H_R(d\bm w)$, and the route entropy to
$\mathcal H_*=-\sum_e\pi^*_e\log\pi^*_e$.

\paragraph{Design.} We use a model that is \emph{only asymptotically}
multivariate regularly varying, so that convergence to the limit is a genuine
test rather than an identity. Route intensities are
$\bm\Lambda=R\,\bm\Theta+\bm\varepsilon$ with $R$ standard Pareto of index
$\tau=1.6$, $\bm\Theta$ drawn from a Dirichlet spectral law $H_R$ on the $\ell_1$
simplex ($d=6$ routes, concentration
$\bm\alpha=(0.6,0.6,0.8,1.0,1.4,2.0)$), and $\bm\varepsilon$ independent
$\mathrm{Exp}(1)$ noise. The light-tailed noise leaves the tail and spectral
measure unchanged but perturbs the direction at finite levels, so the
conditional law approaches the limit only as the conditioning threshold $u$ on
$\|\bm\Lambda\|_1$ grows. Per-route success probabilities are heterogeneous,
$\bm s=(0.15,0.90,0.30,0.60,0.95,0.40)$. Conditional on $\bm\Lambda$, we draw the
first-success route and time from competing exponential clocks with rates
$s_e\Lambda_e$. The limits $\pi^*$, $S^*$, and $\mathcal H_*$ are evaluated by
independent Monte Carlo integration over $H_R$ and the Pareto radius; the model
is simulated with $2\times10^5$ retained extreme draws per threshold.

\begin{figure}[ht]
\centering
\includegraphics[width=\textwidth]{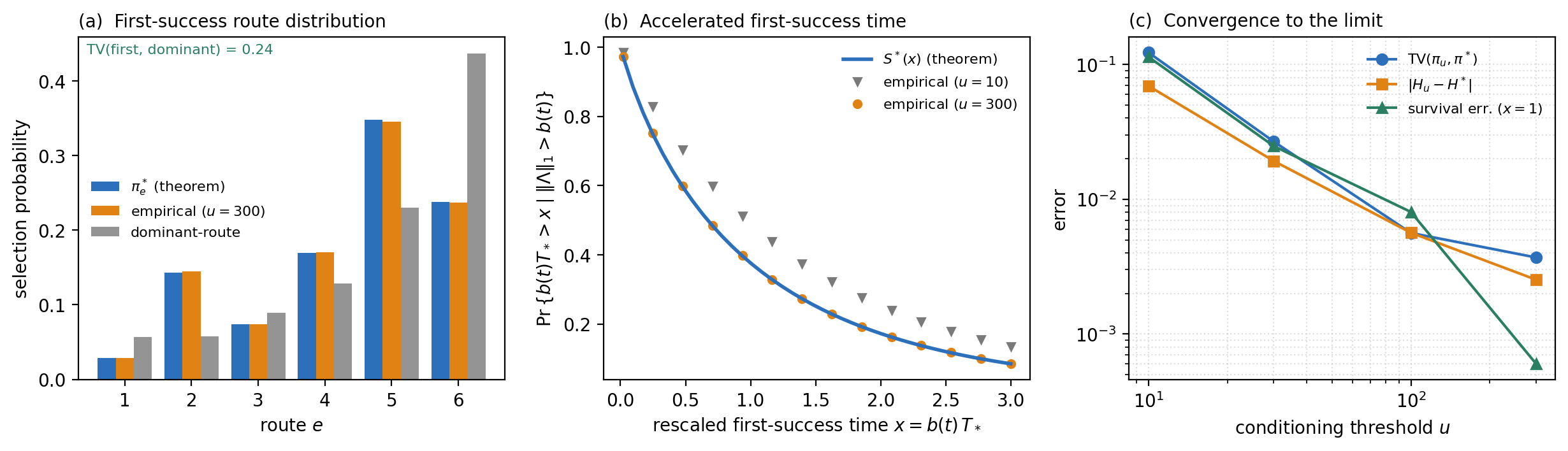}
\caption{Validation of the first successful-event theorem. (a)~First-success
route distribution: theoretical $\pi^*_e$, empirical frequencies at $u=300$, and
the dominant-route distribution, which differs (total variation $0.243$).
(b)~Rescaled first-success time survival: the theoretical curve $S^*(x)$ with
empirical points at $u=10$ (off the limit) and $u=300$ (on the limit).
(c)~Log--log decay of the route total-variation, entropy, and survival errors as
the conditioning threshold $u$ increases.}
\label{fig:clock}
\end{figure}

\paragraph{Results.} Table~\ref{tab:clock} reports convergence across thresholds
$u\in\{10,30,100,300\}$. The total-variation distance between the empirical and
theoretical route distributions falls from $0.123$ to $0.004$, the entropy gap
from $0.069$ to $0.003$, and the first-time survival error at $x=1$ from
$0.114$ to $0.001$. Figure~\ref{fig:clock} shows the agreement at the largest
threshold and the power-law decay of all three error metrics. The first-success
distribution is also numerically distinct from the dominant-route distribution
of the aggregation proposition (total variation $0.243$): the route carrying the
largest flow (route~6, the largest spectral mass) is not the route on which the
first successful event most often occurs (route~5, whose events succeed with
probability $0.95$). This confirms that the success reweighting $s_e$ selects a
genuinely different object.

\begin{table}[ht]
\centering\small
\caption{Convergence of the first successful-event limits. TV is the
total-variation distance between the empirical and theoretical first-success
route distributions; $|\mathcal H_u-\mathcal H_*|$ is the route-entropy gap; the
last columns are absolute errors of the rescaled first-time survival
$\Pr\{b(t)T_*>x\}$ at $x=0.1,0.5,1.0$.}
\label{tab:clock}
\begin{tabular}{rrrrrr}
\toprule
$u$ & TV$(\pi_u,\pi^*)$ & $|\mathcal H_u-\mathcal H_*|$ & err.\ $x{=}0.1$ & err.\ $x{=}0.5$ & err.\ $x{=}1.0$\\
\midrule
10  & $0.1228$ & $0.0694$ & $0.0453$ & $0.1087$ & $0.1140$\\
30  & $0.0268$ & $0.0193$ & $0.0136$ & $0.0266$ & $0.0248$\\
100 & $0.0056$ & $0.0056$ & $0.0034$ & $0.0092$ & $0.0081$\\
300 & $0.0037$ & $0.0025$ & $0.0023$ & $0.0027$ & $0.0006$\\
\bottomrule
\end{tabular}
\end{table}

\section*{S9. Population-tail sensitivity and observational-equivalence twin}

These two synthetic constructions reuse the metapopulation model defined in
main-text Section~4.1 without modification, except for the single parameter that
each construction varies.
Both are reproduced by \texttt{scripts/pooling\_identifiability\_and\_tail.py},
which imports the epidemic generator and SEIR routine of the main simulation
script; a generator check confirms that the reused code is byte-identical to the
primary model at the baseline tail index.

\paragraph{Population-tail sensitivity.} We regenerate the 60-region multiplex with the
population tail index $\alpha$ set to $4.0$, $2.5$, $1.8$, and $1.4$, holding the
seed, gateway, distance, and temporal parameters fixed, and run 120
Full-multiplex SEIR replicates at each value. For each replicate we record its
route composition and total imports, and for each $\alpha$ we compute the
tail-conditioned route entropy and dominant-route share (the spectral
functionals of the main text, conditioning on the upper-quartile outbreaks) and
the pooled-versus-layer-resolved surveillance regret at $k=10$ exactly as in
Section~S7. Table~\ref{tab:tailsweep} reports the sweep. As $\alpha$ decreases,
the generated population distribution becomes more heterogeneous and the
backbone concentrates: the tail-conditioned dominant-route share rises and route
entropy falls monotonically. Changing $\alpha$ regenerates the complete network,
so hub dominance, route allocation, and epidemic dynamics change together. This
finite-network sensitivity does not isolate or validate radial-index
preservation. The Hill estimates describe the realized 60-node, truncated
networks and carry no asymptotic interpretation.
The $\alpha=1.8$ row ($0.394$) reproduces the main-text Full-multiplex
dominant-route share ($0.395$), confirming the sweep harness against the primary
calibration. The
surveillance regret stays between $4$ and $10$ percentage points and is largest
for light tails, because when the tail is heaviest a few dominant hubs lead both
layers and the pooled total-import ranking tracks the air ranking. The
decision cost of pooling therefore persists across these four population-tail
settings.

\begin{table}[ht]
\centering\small
\caption{Population-tail sensitivity of the synthetic epidemic, 120 replicates
per $\alpha$. The Hill column describes the realized 60-node, truncated networks
and does not validate an asymptotic tail law. Route summaries condition on
upper-quartile outbreaks; regret is the air-surveillance gap at $k=10$
(Section~S7). Backbone concentration increases with generated population
heterogeneity, and regret persists across all four settings.}
\label{tab:tailsweep}
\begin{tabular}{rrrrr}
\toprule
$\alpha$ & $\widehat\tau_{\rm Hill}$ & Dominant-route share & Route entropy & Regret at $k{=}10$ (pp)\\
\midrule
$4.0$ & $2.92$ & $0.162$ & $3.62$ & $10.0$\\
$2.5$ & $1.83$ & $0.300$ & $3.46$ & $8.6$\\
$1.8$ & $1.32$ & $0.394$ & $3.06$ & $8.5$\\
$1.4$ & $1.02$ & $0.480$ & $2.67$ & $3.9$\\
\bottomrule
\end{tabular}
\end{table}

\paragraph{Observational-equivalence twin.} The twin holds fixed the effective
pooled coupling $W_{ij}=\kappa_a F^{a}_{ij}+\kappa_c F^{c}_{ij}$ that drives the
epidemic, where $\kappa_a=0.075$ and $\kappa_c=0.012$ are the air and commuting
transmission efficiencies of the main model. Every decomposition of $W$ into air
and commuting layers with the same $W$ generates the identical realized epidemic
and the identical pooled per-destination import counts, so any two such
decompositions are observationally equivalent under pooling. Split~A uses the
model's own layers. Split~B relocates the air share to different destinations by
permuting the columns of the effective air-share field
$\phi_{ij}=\kappa_a F^{a}_{ij}/W_{ij}$, with the permutation chosen from a fixed
search over $2{,}000$ candidates to match Split~A's realized air budget. Holding
the effective coupling fixed, rather than the raw flow, preserves each layer's
per-trip transmissibility, so the air layer keeps its $\approx14\%$ import weight
rather than its $2.6\%$ raw-flow weight. We run 200 replicates of the
shared-coupling SEIR and attribute each destination's realized imports to air
under each split by its $\phi$ weight. The two splits carry a comparable overall
air share ($13.1\%$ and $14.9\%$) but place air on different destinations: their
top-ten air-surveillance targets overlap in five of ten (Jaccard $0.33$), and
the pooled ranking recovers $38.3\%$ of Split~A's air burden at $k=10$ against
$46.6\%$ for the layer-resolved ranking, matching the Section~S7 analysis. This
realizes the noninjectivity of the aggregation map at the epidemic level.

\section*{S10. Proofs of the transfer results}
All vector spaces in this section are finite-dimensional. When the dimension
is $d$, write $\mathbb E_d=[0,\infty)^d\setminus\{\bm0\}$, abbreviated to
$\mathbb E$ when no ambiguity is possible. We use an arbitrary norm and its
induced operator norm; all choices are equivalent in finite dimension. Define
$\mathbb M(\mathbb E_d)$ and $\mathcal C_{\mathbb E_d}$ as in
\eqref{eq:mrv_test_class}. For a measure sequence on $\mathbb E_d$, write
$\mu_t\xrightarrow{M(\mathbb E_d)}\mu$ when
$\int f\,d\mu_t\to\int f\,d\mu$ for every $f\in\mathcal C_{\mathbb E_d}$.
Multivariate regular variation of a random vector $\bm Y$ with index $\tau$
means that, for a scaling function $b\in\mathrm{RV}_{1/\tau}$,
\begin{equation}
t\Prb\{b(t)^{-1}\bm Y\in\cdot\}
\xrightarrow{M(\mathbb E_d)}\mu(\cdot)
\quad\text{on }\mathbb E,
\label{eq:mrv_definition}
\end{equation}
where $\mu$ is a nonzero element of $\mathbb M(\mathbb E_d)$. The limit is
homogeneous: $\mu(rA)=r^{-\tau}\mu(A)$ for $r>0$.

\subsection*{Shared-factor limit and tail coefficient}

\begin{proof}[Proof of Theorem~\ref{thm:shared_factor}]
Let $\nu_\tau$ be the Pareto tail measure on $(0,\infty)$. It satisfies
$\nu_\tau(x,\infty)=x^{-\tau}$ and has density
$\tau r^{-\tau-1}\,dr$. Apply Theorem~A.1 of
\citet{wang2024heterogeneous} with their heavy-tailed variable $X$ equal to
$R$ and their process $\xi(r)$ equal to $\bm A(r)$. The almost-sure
convergence, independence, and uniform $(\tau+\epsilon)$ moment condition are
exactly the hypotheses required there. First,
\begin{equation}
t\Prb\left\{
\left(\bm A(R),\frac{R}{b(t)}\right)\in\cdot
\right\}
\xrightarrow{M}
\Prb(\bm A_\infty\in\cdot)\times\nu_\tau(\cdot)
\label{eq:wr_joint_limit}
\end{equation}
on $[0,\infty)^m\times(0,\infty]$. The moment condition makes the
multiplication map $h(a,r)=ra$ admissible despite its noncompact preimages;
the second part of their theorem therefore gives
\begin{align}
t\Prb\left\{\frac{R\bm A(R)}{b(t)}\in C\right\}
&\longrightarrow
\{\Prb_{\bm A_\infty}\times\nu_\tau\}
\{(a,r):ra\in C\}\nonumber\\
&=
\E\left[\int_0^\infty
\mathbf1\{r\bm A_\infty\in C\}
\tau r^{-\tau-1}\,dr\right]
=\mu_A(C)
\label{eq:wr_product_limit}
\end{align}
for $\mu_A$-continuity sets $C$ bounded away from the origin. This proves
multivariate regular variation of $\bm\Lambda=R\bm A(R)$ and
\eqref{eq:shared_factor_measure}. Lemma~\ref{lem:poissonization}, which
requires conditional Poisson marginals but not conditional independence,
then gives the same limit measure for $\bm Z$.

It remains to calculate the standardized tail coefficient. Positivity and
\eqref{eq:wr_uniform_moment} imply
$0<m_k=\E(A_{\infty,k}^{\tau})<\infty$. From
\eqref{eq:shared_factor_measure}, for $x>0$,
\begin{align}
\mu_A\{y:y_k>x\}
&=\E\int_{x/A_{\infty,k}}^\infty
\tau r^{-\tau-1}\,dr
=x^{-\tau}m_k,
\label{eq:shared_marginal_measure}\\
\mu_A\{y:y_j>xm_j^{1/\tau},\
y_k>xm_k^{1/\tau}\}
&=x^{-\tau}
\E\left[
\min\left\{
\frac{A_{\infty,j}^{\tau}}{m_j},
\frac{A_{\infty,k}^{\tau}}{m_k}
\right\}\right].
\label{eq:shared_joint_measure}
\end{align}
Consequently, each standardized coordinate
$Z_k/m_k^{1/\tau}$ has tail asymptotic
\[
t\Prb\{Z_k>b(t)y\,m_k^{1/\tau}\}
\longrightarrow y^{-\tau},\qquad y>0.
\]
Regular variation of the marginal tails and the standard quantile inversion
property imply that their generalized upper quantiles are asymptotically
$m_k^{1/\tau}b\{(1-q)^{-1}\}$. Dividing the joint limit
\eqref{eq:shared_joint_measure} by either marginal limit and taking
$q\uparrow1$ proves \eqref{eq:shared_factor_lambda}. The same calculation
also shows directly that the expression lies in $[0,1]$, since each
normalized term inside the minimum has expectation one.
\end{proof}

\subsection*{Random nonnegative linear scaling}

The next lemma is a multivariate generalized Breiman result for an independent
random nonnegative matrix. Related matrix-scaling arguments are used for
multivariate degree limits in directed preferential-attachment networks
\citep{wang2022reciprocity}; the process-indexed version used in
Theorem~\ref{thm:shared_factor} is due to
\citet[Appendix~A]{wang2024heterogeneous}.

We use the following random-matrix scaling lemma.
\begin{lemma}
\label{lem:random_scaling}
Suppose $B$ takes values in
$\mathbb E_d=[0,\infty)^d\setminus\{\bm0\}$ and satisfies
\eqref{eq:mrv_definition} with index $\tau$, scaling $b$, and limit $\nu$.
Let $L$ be an independent random nonnegative $m\times d$ matrix with no zero
column almost surely and
$\E\|L\|^{\tau+\epsilon}<\infty$ for some $\epsilon>0$. Then $LB$ is
multivariate regularly varying on
$\mathbb E_m=[0,\infty)^m\setminus\{\bm0\}$ with the same scaling and
\begin{equation}
t\Prb\{b(t)^{-1}LB\in\cdot\}
\xrightarrow{M(\mathbb E_m)}
\mu_L(\cdot),\qquad
\mu_L(A)=\E[\nu\{x:Lx\in A\}].
\label{eq:random_scaling_lemma}
\end{equation}
The measure $\mu_L$ is nonzero, Radon, and homogeneous of order $-\tau$.
\end{lemma}

\begin{proof}
Let $f\in C_c(\mathbb E_m)$. There are constants $0<a<R<\infty$ such that
$f(x)=0$ unless $a\leq\|x\|\leq R$. For every fixed nonnegative matrix
$l$ with no zero column, the function
$x\mapsto f(lx)$ belongs to $C_c(\mathbb E_d)$. To see this using
$\ell_1$ norms, put $s_*(l)=\min_k\sum_r l_{rk}>0$; then
$\|lx\|_1\geq s_*(l)\|x\|_1$ for every $x\geq0$. Hence
\eqref{eq:mrv_definition} gives
\begin{equation}
t\E[f(lB/b(t))]\longrightarrow\int_{\mathbb E_d}f(lx)\,\nu(dx).
\label{eq:fixed_q_limit}
\end{equation}

It remains to justify integration with respect to the law of $L$. Since
$f(lB/b(t))\neq0$ implies $\|B\|>ab(t)/\|l\|$, the radial tail implied by
multivariate regular variation and Potter's bound yield, after increasing a
constant $C$ if necessary,
\begin{equation}
t\E|f(lB/b(t))|
\leq C\|f\|_\infty\{1+\|l\|^{\tau+\epsilon}\}
\label{eq:potter_domination}
\end{equation}
for all sufficiently large $t$ and all such matrices $l$. For completeness,
apply Potter's bound to
$\Prb(\|B\|>b(t)y)/\Prb(\|B\|>b(t))$ with $y=a/\|l\|$ whenever $b(t)y$
exceeds a fixed tail threshold. Otherwise $\|l\|$ is at least of order
$b(t)$, and the right side of \eqref{eq:potter_domination} dominates the
trivial bound $t\|f\|_\infty$ because
$b(t)^{\tau+\epsilon}/t\to\infty$.
The right side of \eqref{eq:potter_domination} is integrable by assumption.
Independence, \eqref{eq:fixed_q_limit}, and dominated convergence therefore
give
\[
t\E[f(LB/b(t))]\longrightarrow
\E\left[\int f(Lx)\,\nu(dx)\right]
=\int f(y)\,\mu_L(dy),
\]
which is $M$-convergence on $\mathbb E_m$.

For any compact $K\subset\mathbb E_m$ with
$\inf_{y\in K}\|y\|=a_K>0$,
\[
\mu_L(K)\leq
\E\nu\{x:\|x\|\geq a_K/\|L\|\}
=C_K\E\|L\|^\tau<\infty,
\]
using homogeneity of $\nu$; thus $\mu_L$ is Radon. Nonnegativity and the
absence of zero columns imply $Lx\neq0$ for every $x\in\mathbb E_d$. If
$A_m=\{y:m^{-1}\leq\|y\|\leq m\}$, then
$\bigcup_{m\geq1}L^{-1}A_m=\mathbb E_d$ almost surely. Since $\nu$ is nonzero,
monotone convergence gives
$\mu_L(\bigcup_m A_m)=\E\nu(\mathbb E_d)>0$, so $\mu_L\neq0$. Finally,
homogeneity follows directly:
$\mu_L(rA)=\E\nu(L^{-1}rA)=r^{-\tau}\mu_L(A)$.
\end{proof}

\subsection*{Poissonization preserves the multivariate tail measure}

The next lemma requires only conditional Poisson marginals, not conditional
independence of the coordinates.

The following lemma shows that conditional Poisson sampling preserves the
multivariate tail measure.
\begin{lemma}
\label{lem:poissonization}
Let $\Lambda$ be multivariate regularly varying on $\mathbb E$ with index
$\tau$, scaling $b$, and limit measure $\mu$. Conditional on $\Lambda$, suppose
that each coordinate $Z_k$ has marginal distribution
$\operatorname{Poisson}(\Lambda_k)$. Then $Z$ is multivariate regularly
varying with the same scaling and the same limit measure $\mu$.
\end{lemma}

\begin{proof}
All norms being equivalent, there is a constant $c_d$ such that
$\|x\|\leq c_d\max_k|x_k|$. For a scalar Poisson variable $P_\lambda$, the
Bernstein-Chernoff inequality gives
\begin{equation}
\Prb(|P_\lambda-\lambda|>u)
\leq2\exp\left\{-\frac{u^2}{2(\lambda+u)}\right\}.
\label{eq:poisson_bernstein}
\end{equation}
Fix $M,\delta>0$. On $\{\|\Lambda\|\leq Mb(t)\}$, a union bound over the
$d$ coordinates and \eqref{eq:poisson_bernstein}, with
$u=\delta b(t)/c_d$, give constants $C_1,C_2>0$ such that
\begin{equation}
\Prb\{\|Z-\Lambda\|>\delta b(t),\ \|\Lambda\|\leq Mb(t)\}
\leq C_1\exp\{-C_2b(t)\}.
\label{eq:exponential_equivalence}
\end{equation}
Because $b\in\mathrm{RV}_{1/\tau}$,
$t\exp\{-C_2b(t)\}\to0$.

Let $f\in C_c(\mathbb E)$ and let $K$ be its support. For sufficiently small
$\delta$, the closed $\delta$-enlargement $K^\delta$ is compact and bounded
away from the origin. On $\{\|Z-\Lambda\|\leq\delta b(t)\}$, the difference
$|f(Z/b(t))-f(\Lambda/b(t))|$ is bounded by the modulus of continuity
$\omega_f(\delta)$ and vanishes unless
$\Lambda/b(t)\in K^\delta$. Consequently,
\begin{align*}
&t\E|f(Z/b(t))-f(\Lambda/b(t))|\\
&\quad\leq
t\omega_f(\delta)\Prb\{\Lambda/b(t)\in K^\delta\}
+2t\|f\|_\infty C_1e^{-C_2b(t)}
+2t\|f\|_\infty\Prb\{\|\Lambda\|>Mb(t)\}.
\end{align*}
Choose $K^\delta$ to be a $\mu$-continuity set, which is possible for all but
countably many enlargement radii. Taking $t\to\infty$, then $M\to\infty$, and
finally $\delta\downarrow0$, the first term vanishes because
$\omega_f(\delta)\to0$, the second by
\eqref{eq:exponential_equivalence}, and the third because
$\mu\{x:\|x\|>M\}=O(M^{-\tau})$. Hence
\[
t\E f(Z/b(t))-t\E f(\Lambda/b(t))\longrightarrow0.
\]
The $M$-convergence for $\Lambda$ proves the same convergence for $Z$.
\end{proof}

\begin{proof}[Proof of Theorem~\ref{thm:structural}]
Apply Lemma~\ref{lem:random_scaling} with $B=\bm B_i$ and $L=Q_i$.
The intensity vector $\Lambda_i=Q_i\bm B_i$ is multivariate regularly varying
with index $\tau$ and limit measure
$\mu_Z(A)=\E[\nu\{b:Q_ib\in A\}]$. Lemma~\ref{lem:poissonization} then shows
that the mixed-Poisson vector $\bm Z_i$ has the same limit measure. The
continuity-set formulation in \eqref{eq:random_pushforward} follows from the
Portmanteau theorem for $M$-convergence.
\end{proof}

\subsection*{Route aggregation, dominant directions, and Gamma mixing}

\begin{proof}[Proof of Proposition~\ref{prop:routes}]
Every column of $P$ sums to one. Because $Q^{R}$ has positive diagonal
entries, the nonnegative matrix $Q^{R}P$ has no zero column almost surely.
Moreover,
$\|Q^{R}P\|\leq\|Q^{R}\|\,\|P\|$, so its
$(\tau+\epsilon)$th moment is finite. Lemma~\ref{lem:random_scaling}
therefore shows that the route-intensity vector
$Q^{R}P\bm B^{A}$ is multivariate regularly varying with limit
\eqref{eq:route_limit}. Lemma~\ref{lem:poissonization} then gives
\[
t\Prb\{\bm Z^{R}/b(t)\in\cdot\}
\xrightarrow{M(\mathbb E_{|\mathcal R|})}\mu_R(\cdot).
\]
We first prove the aggregation statement. Equip both finite-dimensional
spaces with the $\ell_1$ norm. Because $S$ is nonnegative and has no zero
column,
\[
\|Sz\|_1=\sum_{k}\left(\sum_r S_{rk}\right)z_k
\geq s_*\|z\|_1,\qquad
s_*=\min_k\sum_rS_{rk}>0
\]
for $z\in[0,\infty)^{|\mathcal R|}$. Hence the preimage under $S$ of
any compact set bounded away from zero is bounded and bounded away from zero.
The map $z\mapsto Sz$ is continuous, so the mapping theorem for
$M$-convergence gives
\[
t\Prb\{S\bm Z^{R}/b(t)\in\cdot\}
\xrightarrow{M}\mu_R\circ S^{-1}(\cdot)=\mu_S(\cdot).
\]
The same lower bound shows that $\mu_S$ is Radon. Its homogeneity follows
from the homogeneity of $\mu_R$ and linearity of $S$. Moreover,
$Sz\ne0$ for every nonzero $z\geq0$, so the pushforward of the nonzero
measure $\mu_R$ is nonzero. It is therefore a valid
multivariate-regular-variation limit.

For the second statement, let
$\mathbb S_+=\{w\geq0:\|w\|_1=1\}$ and define
\[
H_R(A)=
\frac{\mu_R\{z:\|z\|_1>1,
z/\|z\|_1\in A\}}
{\mu_R\{z:\|z\|_1>1\}},\qquad A\subset\mathbb S_+.
\]
The denominator is finite and positive. Indeed, finiteness follows from the
radial form of a homogeneous Radon measure, and positivity follows from
nonzero homogeneity of $\mu_R$. Polar decomposition of a homogeneous Radon
measure gives
\begin{equation}
\Prb\!\left\{\bm Z^{R}/\|\bm Z^{R}\|_1\in A
\,\middle|\, \|\bm Z^{R}\|_1>u\right\}
\longrightarrow H_R(A)
\label{eq:spectral_convergence_proof}
\end{equation}
for every $H_R$-continuity set $A$ as $u\to\infty$.

Let $T$ be the subset of $\mathbb S_+$ on which two or more coordinates
tie for the maximum. By assumption $H_R(T)=0$. Each $C_e$ is therefore an
$H_R$-continuity set: its relative boundary is contained in $T$. Outside
$T$, the direction belongs to exactly one $C_e$, and the largest coordinate
of the direction is the largest coordinate of the original vector. Applying
\eqref{eq:spectral_convergence_proof} to $C_e$ proves
\eqref{eq:dominant_route}; the event on which the tie-breaking rule matters
has directions in $T$ and conditional probability tending to zero.

Write $\pi_e(u)$ for the conditional probability on the left of
\eqref{eq:dominant_route} and $\pi_e=H_R(C_e)$. The set
$\mathcal R$ is finite, so coordinatewise convergence
$\pi_e(u)\to\pi_e$, together with continuity of
$x\mapsto-x\log x$ on $[0,1]$ under the convention $0\log0=0$,
gives
$-\sum_e\pi_e(u)\log\pi_e(u)\to\mathcal H_R$.
\end{proof}

\begin{proof}[Proof of Corollary~\ref{cor:aggregation_spectral}]
Write the polar decomposition of $\mu_R$ as
\[
\mu_R(dz)
=c_R\tau r^{-\tau-1}\,dr\,H_R(dw),
\qquad z=rw,\quad r>0,\quad w\in\mathbb S_+,
\]
where $c_R=\mu_R\{z:\|z\|_1>1\}\in(0,\infty)$. For a Borel set
$A$ on the aggregate unit simplex whose boundary has zero transformed
spectral mass,
\begin{align*}
&\mu_S\left\{y:\|y\|_1>1,\
\frac{y}{\|y\|_1}\in A\right\}\\
&\quad=
c_R\int_{\mathbb S_+}\int_0^\infty
\mathbf1\{r\|Sw\|_1>1\}
\mathbf1\{T_S(w)\in A\}
\tau r^{-\tau-1}\,dr\,H_R(dw)\\
&\quad=
c_R\int_{\mathbb S_+}
\|Sw\|_1^\tau
\mathbf1\{T_S(w)\in A\}\,H_R(dw).
\end{align*}
Taking $A$ to be the entire aggregate simplex gives the normalizing
constant. Division proves \eqref{eq:aggregate_spectral}.

For nonidentifiability, suppose $w_1\ne w_2$ and
$T_S(w_1)=T_S(w_2)$. The two route spectral measures
$H_R^{(1)}=\delta_{w_1}$ and $H_R^{(2)}=\delta_{w_2}$ are distinct. In
\eqref{eq:aggregate_spectral}, the factor $\|Sw_i\|_1^\tau$ appears in both
the numerator and denominator and cancels. Both measures therefore induce
the same point mass $\delta_{T_S(w_1)}$ after aggregation. Thus $H_S$ cannot
identify $H_R$ without an additional restriction that separates the fibers
of $T_S$.
\end{proof}

\begin{proof}[Proof of Proposition~\ref{prop:identifiability}]
Because $S$ is nonnegative with no zero column, $\|s_e\|_1>0$ for every $e$, and
for $w\in\mathbb S_+$,
\[
Sw=\sum_e w_e s_e=\sum_e w_e\|s_e\|_1\,\hat s_e,\qquad
\|Sw\|_1=\mathbf 1^\top Sw=\sum_e w_e\|s_e\|_1,
\]
the second equality using $Sw\ge0$ and $\mathbf 1^\top s_e=\|s_e\|_1$. Define the
reweighting map $\Phi(w)_e=w_e\|s_e\|_1/\sum_f w_f\|s_f\|_1$. Since all
$\|s_e\|_1>0$, $\Phi$ is a bijection of $\mathbb S_+$ onto itself, with inverse
$\Phi^{-1}(u)_e=u_e\|s_e\|_1^{-1}/\sum_f u_f\|s_f\|_1^{-1}$. Substituting,
\begin{equation}
T_S(w)=\frac{Sw}{\|Sw\|_1}=\sum_e\Phi(w)_e\,\hat s_e=\hat S\,\Phi(w),
\qquad \hat S=[\hat s_1,\dots,\hat s_{|\mathcal R|}].
\label{eq:TS_factorization}
\end{equation}
As $\Phi$ is a bijection, $T_S$ is injective on $\mathbb S_+$ if and only if the
affine map $u\mapsto\hat S u=\sum_e u_e\hat s_e$ is injective on $\mathbb S_+$.
The latter fails exactly when there are coefficients $c\ne0$ with
$\sum_e c_e\hat s_e=0$ and $\sum_e c_e=0$: given such $c$, its positive and
negative parts, rescaled to unit mass, are distinct simplex points with the same
image; conversely two simplex points with the same image give such a $c$ by
subtraction. Absence of such $c$ is the definition of affine independence of
$\{\hat s_e\}$. This proves the equivalence.

For part~(i), the $\hat s_e$ lie in the affine hyperplane
$\{x:\mathbf 1^\top x=1\}$ of $\mathbb R^m$, which has dimension $m-1$; an
affinely independent subset of an $(m-1)$-dimensional affine space has at most
$m$ points, so injectivity forces $|\mathcal R|\le m$, and $m=1$ forces
$|\mathcal R|=1$. For part~(ii), injectivity makes $T_S$ a homeomorphism onto its
image with inverse $\psi$. Taking $f(v)=\varphi(\psi(v))\|S\psi(v)\|_1^{-\tau}$
in \eqref{eq:aggregate_spectral} and using $\psi\circ T_S=\mathrm{id}$ gives
$\int f\,dH_S=Z^{-1}\int\varphi\,dH_R$ with
$Z=\int_{\mathbb S_+}\|Sw\|_1^\tau H_R(dw)$; the choice $\varphi\equiv1$
evaluates $Z^{-1}$ as the denominator of \eqref{eq:spectral_recovery}, proving
recovery.

For the layer-pooling statement, $s_{(i,j,\ell)}=e_{(i,j)}$ for every route, so
each column of $S$ has $\ell_1$ norm one. Hence, for every $w\in\mathbb S_+$,
\[
\|Sw\|_1=\sum_{i,j,\ell}w_{ij\ell}=1,
\qquad
(Sw)_{ij}=\sum_\ell w_{ij\ell}.
\]
The weight $\|Sw\|_1^\tau$ in \eqref{eq:aggregate_spectral} is therefore one,
and that equation reduces to
\[
H_S(A)=H_R\{w:Sw\in A\}=H_{\rm OD}(A),
\]
which proves \eqref{eq:od_pushforward}. If two routes share $(i,j)$ and differ
only in $\ell$, their signatures coincide. The full map $T_S$ is then
noninjective, so the within-pair layer composition is not identified even though
the source-destination pushforward is.
\end{proof}

\begin{proof}[Proof of Theorem~\ref{thm:first_invasion}]
Put
\[
R_t=\frac{\|\bm\Lambda\|_1}{b(t)},\qquad
W_t=\frac{\bm\Lambda}{\|\bm\Lambda\|_1}.
\]
Polar decomposition of the regularly varying limit measure gives the
conditional convergence
\begin{equation}
\mathcal L\{(R_t,W_t)\mid\|\bm\Lambda\|_1>b(t)\}
\Longrightarrow
\nu_\tau^{(1)}\times H_R,
\label{eq:conditional_polar_limit}
\end{equation}
where $\nu_\tau^{(1)}$ is the Pareto probability distribution on
$[1,\infty)$ with density $\tau r^{-\tau-1}$. In particular, the limiting
radius and angle are independent.

Conditional on $\bm\Lambda$, thinning infectious arrivals by their
establishment probabilities produces independent Poisson processes with
rates $s_e\Lambda_e$. The competing-exponential formula therefore gives
\begin{align}
\Prb(E_*=e\mid\bm\Lambda)
&=\frac{s_e\Lambda_e}{s^\top\bm\Lambda},
\label{eq:competing_route}\\
\Prb(E_*=e,T_*>y\mid\bm\Lambda)
&=\frac{s_e\Lambda_e}{s^\top\bm\Lambda}
\exp\{-y\,s^\top\bm\Lambda\},
\label{eq:competing_route_time}\\
\Prb(T_*>y\mid\bm\Lambda)
&=\exp\{-y\,s^\top\bm\Lambda\}.
\label{eq:competing_time}
\end{align}
The functions
\[
w\longmapsto\frac{s_ew_e}{s^\top w}
\quad\text{and}\quad
(r,w)\longmapsto
\frac{s_ew_e}{s^\top w}e^{-xr s^\top w}
\]
are bounded and continuous except possibly on $\{w:s^\top w=0\}$, which
has $H_R$-measure zero by \eqref{eq:positive_establishment_mass}. Applying
\eqref{eq:conditional_polar_limit} to \eqref{eq:competing_route} proves
\eqref{eq:first_route_limit}. In
\eqref{eq:competing_route_time}, set $y=x/b(t)$ and note that
\[
y\,s^\top\bm\Lambda=xR_t\,s^\top W_t.
\]
Another application of \eqref{eq:conditional_polar_limit} gives
\eqref{eq:first_route_time_limit}. Summing over the finite route set, or
applying the same argument directly to \eqref{eq:competing_time}, proves
\eqref{eq:first_time_limit}.

Let $\pi_e(t)$ denote the conditional first-route probability on the
left of \eqref{eq:first_route_limit}. Coordinatewise convergence
$\pi_e(t)\to\pi_e^*$ and finiteness of $\mathcal R$, together with continuity
of $-x\log x$ on $[0,1]$ under $0\log0=0$, yield
\[
-\sum_e\pi_e(t)\log\pi_e(t)
\longrightarrow-\sum_e\pi_e^*\log\pi_e^*=\mathcal H_*.
\]
\end{proof}

\begin{proof}[Proof of Corollary~\ref{cor:nb}]
Conditional on $(\bm B_*,M,G)$, the count vector has Poisson marginals with
mean vector $GM\bm B_*$. The matrix $GM$ is nonnegative, has no zero column,
is independent of $\bm B_*$, and satisfies the moment condition of
Lemma~\ref{lem:random_scaling}. That lemma followed by
Lemma~\ref{lem:poissonization} proves multivariate regular variation with
limit \eqref{eq:nb_pushforward}. If
$G_k\sim\operatorname{Gamma}(\kappa_k,\kappa_k)$, integrating a
Poisson distribution with mean $G_k(M\bm B_*)_k$ over $G_k$ gives the stated
negative-binomial marginal.

If $G=G_0I$, first condition on $M$. Homogeneity of $\nu_*$ gives
\begin{align*}
\mu_{GM}(A)
&=\E\left[\nu_*\{b:G_0Mb\in A\}\right]\\
&=\E\left[G_0^\tau
\nu_*\{b:Mb\in A\}\right]\\
&=\E(G_0^\tau)\mu_M(A),
\end{align*}
where the last equality uses independence of $G_0$ and $M$. Multiplication of
a homogeneous limit measure by a positive constant cancels in its spectral
normalization, proving the common-multiplier claim. No such scalar
factorization holds in general for coordinate-specific multipliers.
\end{proof}

\bibliographystyle{plainnat}
\bibliography{references}